\documentclass[12pt]{article}
\usepackage{amsmath, amssymb}
\begin{document}
\centerline{\bf Techniques for the Analytic Proof of}\centerline{\bf the Finite Generation of the Canonical Ring}

\bigbreak\centerline{Yum-Tong Siu\ %
\footnote{Partially supported by Grant 0500964 of the National Science
Foundation.} }

\bigbreak\centerline{\it Written for the Proceedings of}\centerline{\it the Conference on Current Developments in Mathematics in}\centerline{\it Harvard University, November 16-17, 2007}

\bigbreak

\bigbreak\noindent{\sc Table of Contents}
\begin{itemize}\item[\S1.] Introduction
\item[\S2.] Reduction of Finite Generation of Canonical Ring to Achievement of Stable Vanishing Order
\item[\S3.]  Decomposition of Closed Positive (1,1)-Currents and their Modified Restrictions to Hypersurfaces 
\item[\S4.] Discrepancy Subspaces
\item[\S5.] Construction of Pluricanonical Sections with Fixed Sufficiently Ample Twisting
\item[\S6.] Subspaces of Minimum Additional Vanishing for the Second Case of the Dichotomy
\item[\S7.] Big Sum of a Line Bundle and the Canonical Line Bundle
\end{itemize}

\bigbreak

\bigbreak\noindent{\sc \S1. Introduction.}  This article is an exposition of the analytic proof of the finite generation of the canonical ring for a compact complex algebraic manifold of general type [Siu 2006, Siu 2007, Siu 2008].  An algebraic proof was given in [Birkan-Cascini-Hacon-McKernan 2006].

\bigbreak\noindent(1.1) {\it Main Theorem.} Let $X$ be a compact complex
algebraic manifold of general type.  Then the canonical ring
$$
R\left(X, K_X\right)=\bigoplus_{m=1}^\infty\Gamma\left(X,
mK_X\right)
$$
is finitely generated, where $K_X$ is the canonical line bundle of $X$.

\bigbreak In this exposition we will list and discuss the main techniques and explain how they are put together in the proof.   Of the various main techniques some special attention is given to
\begin{itemize}\item[(i)] the technique of discrepancy subspaces and
\item[(ii)] the technique of subspaces of ``minimum additional vanishing''
\end{itemize}
respectively treated in \S4 and \S6.

\medbreak The technique of discrepancy subspaces was already discussed in detail in [Siu 2008].  It measures the deviation from a sufficiently ample line bundle and is used to terminate the process of possibly infinitely many blow-ups.  Here we explain its motivation from the perspective of dynamic multiplier ideal sheaves. When multiplier ideal sheaves were first introduced by Kohn [Kohn 1979] as measurements of failure of estimates in partial differential equations and introduced by Nadel [Nadel 1990] as destabilizing sheaves, their definitions are formulated from the most crucial estimates and involve respectively a family of inequalities and a sequence of inequalities.  They are {\it dynamic} in contrast to the usual multiplier ideal sheaves used in algebraic geometry which translate Nadel's destabilizing subsheaves when the sequence of inequalities used in the definition becomes a single one.  We detail here in \S4 how the notion of discrepancy subspaces arises from the most crucial estimates for the analytic proof of the finite generation of the canonical ring.

\medbreak The technique of subspaces of ``minimum additional vanishing'', treated in \S6,  is used to handle the extension of sections from the second case of the dichotomy (see (3.2)) and has so far been discussed only with very sparse details elsewhere.

\medbreak At the end of this article there is a very brief discussion of the modification needed for the analytic proof of the twisted case of the finite generation of the canonical ring.

\bigbreak\noindent(1.2) {\it List of Main Techniques.}  Here is a list of the main techniques used in the analytic proof of the finite generation of the canonical ring, with numbering from (A) to (H).

\bigbreak\noindent(A) The metric $e^{-\varphi}=\frac{1}{\Phi}$ of minimum singularity for the canonical line bundle $K_X$, which is constructed from an infinite sum $\Phi$ of the $m$-th root of the absolute-value-squares of  the ${\mathbb C}$-basis elements of $\Gamma\left(X, mK_X\right)$ for all $m\in{\mathbb N}$ (see (2.5)). 

\bigbreak\noindent(B) The application of Skoda's theorem for ideal generation to a Zariski open subset of $X$ which can be regarded as a Stein domain spread over ${\mathbb C}^n$ so that the finite generation of the canonical ring is reduced to the precise achievement of all the stable vanishing orders, which means that the infinite sum $\Phi$ is comparable to one of its finite partial sum in the sense that each one is bounded by a positive constant times the other (see (2.6)).   

\medbreak A stable vanishing order means the vanishing order of $\Phi$ across some subvariety at some generic point and in a generic normal direction, which more precisely means the vanishing order of $\Phi$ across some nonsingular hypersurface at a generic point of the hypersurface when the subvariety is blown up to give nonsingular hypersurfaces in normal crossing and $\Phi$ is an infinite sum formed for the blown-up manifold.

\bigbreak\noindent(C) The unique decomposition of a closed positive $(1,1)$-current on an open subset $U$ of ${\mathbb C}^n$ as the sum of a possibly infinite sum of distinct irreducible hypersufaces with positive real coefficients and a residue current whose Lelong number is zero outside at most a countable union of subvarieties of complex codimension at least two in $U$ (see (3.2)).   

\medbreak There is a dichotomy into two cases.  The first case is when either the ${\mathbb R}$-linear combination of distinct irreducible hypersurfaces contains an infinite number of terms or the residue current is not identically zero.  The second case is when there are only a finite number of terms in  the ${\mathbb R}$-linear combination of distinct irreducible hypersurfaces and the residue current is identically zero.

\bigbreak\noindent(D)  When some stable vanishing order is assumed to be not achievable, the theorem of Kawamata-Viehweg-Nadel [Kawamata 1982, Viehweg 1982, Nadel 1990] for multiplier ideal sheaves is applied to the subspace of ``mininum additional vanishing'' (after being blown up to a nonsingular hypersurface) to extend a section of the modified pluricanonical line bundle on it to $X$ (see \S6).  Such an extension would give the precise achievement of the stable vanishing order.  The modified pluricanonical line bundle means the canonical line bundle minus the hypersurface times the stable vanishing order.  The section to be extended may need to be constructed from the second case of the dichotomy in the decomposition of the modified curvature current of the canonical line bundle on the subspace of ``mininum additional vanishing'' (which is now a hypersurface after the blow-up).  

\medbreak This kind of subspace of ``mininum additional vanishing'' is the analog of the minimum center of log canonical singularity [Kawamata 1985, Shokurov 1985] in the techniques used in the study of the Fujita conjecture [Fujita 1987] and in Shokurov's non-vanishing theorem [Shokurov 1985].  Unlike the situations in the Fujita conjecture and in Shokurov's non-vanishing theorem, in the analytic proof of the finite generation of the canonical ring there is already the stable vanishing order and the subspace is to be defined by more vanishing order but only minimally more.   

\medbreak One technique of applying the vanishing theorem of Kawamata-Viehweg-Nadel for multiplier ideal sheaves is applied to certain subspaces $S$ of $X$ is to apply it indirectly through $X$.  We use the vanishing of two sandwiching cohomology groups on $X$ in the long cohomology exact sequence from the short exact sequence for the ideal sheaf of $S$ to get the vanishing of cohomology on $S$.  The reason is that though $X$ is of general type, there is no information about the canonical line bundle of the subspace $S$.

\medbreak Some positive lower bound for the curvature current can be obtained from the general type property of $X$, because a sufficiently small modification of the construction of the metric in question would not affect its multiplier ideal sheaf.  Moreover, the canonical line bundle of the ambient space $X$ rather than that of the subspace $S$ should be in the sheaf used in the cohomology group.   

\bigbreak\noindent(E) The technique of constructing positively curved metrics with additional high-order singularity on a hypersurface of the first case of the dichotomy (see (5.5)).  The construction uses the contribution from the round-up in the computation of the multiplier ideal sheaf and uses a sufficiently ample line bundle to guarantee a positive growth order for the dimension of the space of sections of amply twisted multiples of the modified pluricanonical line bundle.   Here Kronecker's theorem on diophantine approximation is used (see (5.2)).  The general type property of $X$ is used to take care of the ample twisting after taking a high-order roots of the sections of the amply twisted multiples of the modified pluricanonical line bundle.  These metrics are used in extending sections on the subvariety of additional high-order singularity to the ambient manifold $X$.

\bigbreak\noindent(F)  Shokurov's technique of comparing two results from the theorem of Hirzebruch-Riemann-Roch, one applied to a line bundle and the other applied to the line bundle twisted by a flat line bundle.  Shokurov originally introduced his technique for the proof of his non-vanishing theorem [Shokurov 1985]

\medbreak Here it is applied to a subspace of ``minimum additional vanishing'' (see (6.4)).  The flat line bundle occurs, because the second case of the dichotomy involves only the curvature current which determines the line bundle only up to an additive flat line bundle.

\bigbreak\noindent(G) Kronecker's theorem on diophantine approximation. Besides its use in the construction of metrics in (5.2) as described in (E), it is also used to show that the stable vanishing orders are rational, otherwise the round-up and round-down process in the construction of multiplier ideal sheaves would result in the decrease of some stable vanishing orders, which is not possible.

\bigbreak\noindent(H) The notion of a discrepancy subspace, which measures the extent of
failure of achieving stable vanishing order with appropriately-defined uniformity
in $m$ for all $m$-canonical bundles (see (4.2)).  This uniformity is obtained by measuring the deviation from a line bundle which is ample enough for the global generation of multiplier ideal sheaves.  
Its formulation is motivated by the original idea of defining multiplier ideal sheaves from the crucial estimates.

\medbreak Discrepancy subspaces are introduced to handle the problem of infinite number of interminable blow-ups in the process of proving the precise achievement of stable vanishing order.  The termination of the infinite process of blow-ups in the proof is essential, because the extension of section can only be done from a subspace of ``minimum additional singularity'' after we blow up the subspace.    

\medbreak This technique of measuring the deviation from a sufficiently ample line bundle was already introduced in the proof of the deformational invariance of plurigenera [Siu 1998, Siu 2002], though it was not exactly so described in the proof of the deformational invariance of plurigenera.   The independence on $m$ allows also the division of the fixed sufficiently ample line bundle by $m$ so that its contribution is so minimized that its removal makes no difference in the construction of the multiplier ideal sheaf.

\medbreak The method of pluricanonical extension (which is obtained by this technique) used in the proof of the deformational invariance of plurigenera holds the key to the proof, either analytic or algebraic-geometric, of the finite generation of the canonical ring in making it possible to implement some process on induction on the dimension by restricting a pluricanonical section to the base-point set.   
Here pluricanonical extension itself is not explicitly used.  Instead the key ingredients of its proof are directly used here.    

\medbreak When we inductively use discrepancy subspaces of lower dimension, we have to introduce some form of ``holomorphic fibration'' with fibers defined by multiplier ideal sheaves of metrics which are defined by multi-valued holomorphic sections vanishing to high order at a generic point of the fiber.  Multiplier ideal sheaves are used to define such a fibration because they are defined in such a way that they are unchanged by slight perturbations of the metrics  (see (5.5.1), (5.5.2), and (5.6)).  

\bigbreak\noindent(1.3) {\it Notations.} The notations ${\mathbb C}$, ${\mathbb R}$, ${\mathbb Q}$, and ${\mathbb N}$ denote respectively the complex numbers, the real numbers, the rational numbers, and the positive integers.  

\medbreak For a real number $u$ the expression $\left\lfloor u\right\rfloor$ means the round-down of $u$ which is the largest integer not exceeding $u$.   The expression $\left\lceil u\right\rceil$ means the round-up of $u$ which is the smallest integer not smaller than $u$.  

\medbreak  The reduced structure sheaf of a complex space $W$ is denoted by
${\mathcal O}_W$.  The stalk of ${\mathcal O}_W$ at a point $P$ of $W$ is denoted by ${\mathcal O}_{W,P}$.  The maximum ideal of a point $P$ of $W$
is denoted by ${\mathfrak m}_{W,P}$.  The canonical line bundle of a complex manifold $Y$ is denoted by $K_Y$.  The full ideal sheaf of a
subvariety $Z$ in a complex manifold $Y$ is denoted by ${\mathcal
I}d_Z$.  The canonical section of the line bundle associated to a
complex hypersurface $V$ in a complex manifold $Y$ is denoted by
$s_V$ (whose divisor is $V$).

\medbreak A multi-valued
holomorphic section $\sigma$ of a holomorphic line bundle $E$ over a
complex manifold $Y$ means that $\sigma^m$ is a holomorphic section
of $mE$ over $Y$ for some positive integer $m$.    The statement that the vanishing order of $\sigma$
at a point $P_0$ (respectively across a regular hypersurface $Z$) is $q$ means that the vanishing order of $\sigma^m$
at a point $P_0$ (respectively across a regular hypersurface $Z$) is $mq$.   

\medbreak When a ${\mathbb Q}$-divisor $\alpha Y$, with $\alpha\in{\mathbb Q}$ and $Y$ being an irreducible hypersurface, is multiplied by $p\in{\mathbb N}$ to become a holomorphic line bundle, the integer $p$ is automatically chosen with the property that $p\alpha$ is an integer and this choice of $p$ will not be explicitly mentioned and is understood.

\bigbreak\noindent{\sc \S2. Reduction of Finite Generation of Canonical Ring to Achievement of Stable Vanishing Order.}  

\bigbreak\noindent(2.1) {\it Definition of Multiplier Ideal
Sheaves.} For a local plurisubharmonic function $\varphi$ on an open
subset of ${\mathbb C}^n$, the multiplier ideal sheaf ${\mathcal
I}_\varphi$ is the sheaf of germs of holomorphic functions $f$ such
that $\left|f\right|^2e^{-\varphi}$ is locally integrable.

\medbreak\noindent(2.1.1) {\it Remark.} This is the usual definition of multiplier ideal sheaf in algebraic geometry which is defined by a single inequality and is {\it static}, as opposed to a {\it dynamic} multiplier ideal sheaf which is defined by a family of inequalities.

\bigbreak\noindent(2.2) {\it Statement on Global Generation of
Multiplier Ideal Sheaves} [Siu 1998]. Let $L$ be a holomorphic line bundle over
an $n$-dimensional compact complex manifold $Y$ with a Hermitian
metric which is locally of the form $e^{-\xi}$ with $\xi$
plurisubharmonic. Let ${\cal I}_\xi$ be the multiplier ideal sheaf
of the Hermitian metric $e^{-\xi}$. Let $A$ be an ample holomorphic
line bundle over $Y$ such that for every point $P$ of $Y$ there are
a finite number of elements of $\Gamma(Y,A)$ which all vanish to
order at least $n+1$ at $P$ and which do not simultaneously vanish
outside $P$. Then $\Gamma(Y,{\cal I}_\xi\otimes(L+A+K_Y))$ generates
${\cal I}_\xi\otimes(L+A+K_Y)$ at every point of $Y$.

\bigbreak\noindent(2.3) {\it Skoda's Result on Ideal Generation}
[Skoda 1972]. Let $\Omega $ be a domain spread over $\mathbb{C}^{n}$
which is Stein. Let $\psi$ be a plurisubharmonic function on $\Omega
$, $g_{1},\ldots ,g_{p}$ be holomorphic functions on $\Omega $,
$\alpha
>1$, $q=\min \left(n,p-1\right)$, and $f$ be a holomorphic function on $\Omega
$. Assume that
$$\int_{\Omega }\frac{\left\vert f\right|^{2}e^{-\psi }}{\left(
\sum_{j=1}^p\left| g_{j}\right| ^{2}\right)^{\alpha q+1}}<\infty.$$
Then there exist holomorphic functions $h_{1},\ldots ,h_{p}$ on
$\Omega$ with $f=\sum_{j=1}^{p} h_jg_j$ on $\Omega$ such that
$$\int_{\Omega }\frac{\left\vert h_{k}\right|^{2}e^{-\psi
}}{\left(\sum_{j=1}^p\left|g_{j}\right|^2\right)^{\alpha q}}\leq
\frac{\alpha }{\alpha -1}\int_{\Omega
}\frac{\left|f\right|^2e^{-\psi }}{\left( \sum_{j=1}^p\left|
g_{j}\right|^2\right)^{\alpha q+1}}$$ for $1\leq k\leq p$.

\medbreak\noindent(2.3.1) {\it Remark.}  Skoda's original statement is for a Stein domain $\Omega$ in ${\mathbb C}^n$, but his proof works also for a domain spread over $\mathbb{C}^{n}$
which is Stein.  We need the setting of a domain $\Omega$ spread over $\mathbb{C}^{n}$, which will be a Zariski open subset of our compact complex algebraic manifold $X$ with a finite-to-one holomorphic map $\pi: X\to{\mathbb P}_n$ and will spread over the affine part ${\mathbb C}^n$ of ${\mathbb P}_n$ under $\pi$.

\bigbreak\noindent(2.4) {\it Multiplier-Ideal Version of Skoda's
Result on Ideal Generation.} Let $X$ be a compact complex algebraic
manifold of complex dimension $n$, $L$ be a holomorphic line bundle
over $X$, and $E$ be a holomorphic line bundle on $X$ with metric
$e^{-\psi }$ such that $\psi $ is plurisubharmonic. Let $k\geq 1$ be
an integer, $G_{1},\ldots ,G_{p}\in \Gamma (X,L)$, and
$\left|G\right|^{2}=\underset{j=1}{\overset{p}{\sum
}}\left|G_j\right|^{2}$. Let
$\mathcal{I=I}_{(n+k+1)\log\left|G\right|^{2}+\psi}$ and
$\mathcal{J=I}_{(n+k)\log\left|G\right|^2+\psi}$. Then
$$\displaylines{\qquad\qquad\Gamma\left(X,\mathcal{I}\otimes
\left((n+k+1)L+E+K_{X}\right)\right)\hfill\cr\hfill
=\underset{j=1}{\overset{p}{%
\sum }}G_{j}\,\Gamma \left(X,\mathcal{J}\otimes\left(
(n+k)L+E+K_{X}\right)\right).\qquad\qquad\cr}$$

\medskip\noindent {\it Proof}. Take $F\in\Gamma\left(X,\mathcal{I}\otimes
\left((n+k+1)L+E+K_{X}\right)\right)$.  Let $S$ be a meromorphic
section of $E$. Take a branched cover map $\pi:X\rightarrow{\bf
P}_n$.  Let $Z_0$ be a hypersurface in ${\bf P}_n$ which contains
the infinity hyperplane of ${\bf P}_n$ and the branching locus of
$\pi$ in ${\bf P}_n$ such that $Z:=\pi^{-1}(Z_0)$ contains the
divisor of $G_1$ and both the pole-set and zero-set of $S$ . Let
$\Omega=X-Z$. Let $g_j=\frac{G_j}{G_1}$ ($1\leq j\leq p$) and
$\left|g\right|^2=\sum_{j=1}^p\left|g_j\right|^2$.  Define $f$ by
$$
\frac{F}{G_1^{n+k+1}S}=fdz_1\wedge\cdots\wedge dz_n,
$$
where $z_1,\cdots,z_n$ are the affine coordinates of ${\mathbb
C}^n$.  Use $\alpha=\frac{n+k}{n}$.  Let $\varphi=\psi-\log|S|^2$.
It follows from $F\in{\cal I}_{(n+k+1)\log|G|^2+\psi}$ locally
that
$$
\int_X\frac{|F|^2}{|G|^{2(n+k+1)}}e^{-\psi}<\infty,
$$
which implies that
$$
\int_\Omega\frac{|f|^2}{|g|^{2(n+k+1)}}e^{-\varphi}=
\int_\Omega\frac{\left|\frac{F}{G_1^{n+k+1}S}\right|^2}{\left|\frac{G}{
G_1}\right|^{2(n+k+1)}}e^{-\varphi}=
\int_\Omega\frac{|F|^2}{|G|^{2(n+k+1)}}e^{-\psi}<\infty.
$$
By Skoda's theorem on ideal generation (2.3) with $q=n$ (which we
assume by adding some $F_{p+1}\equiv\cdots\equiv F_{n+1}\equiv 0$ if
$p<n+1$) so that $2\alpha q+2=2\cdot\frac{n+k}{n}\cdot
n+2=2(n+k+1)$, we obtain holomorphic functions $h_1,\cdots,h_p$ on
$\Omega$ such that $f=\sum_{j=1}^p h_jg_j$ and
$$
\sum_{j=1}^p\int_\Omega\frac{|h_j|^2}{|g|^{2(n+k)}}e^{-\varphi}<\infty.
$$
Define
$$
H_j=G_1^{n+k}h_j S dz_1\wedge\cdots\wedge dz_n.
$$
Then $F=\sum_{j=1}^p H_j G_j$ and
$$
\int_\Omega\frac{|H_j|}{|G|^{2(n+k)}}e^{-\psi}=
\int_\Omega\frac{|h_j|}{|g|^{2(n+k)}}e^{-\varphi}<\infty
$$
so that $H_j$ can be extended to an element of
$\Gamma(X,(n+k)L+E+K_X)$. Q.E.D.

\bigbreak\noindent(2.5) {\it Metric of Minimum Singularity and Definition of Precise Achievement of Stable Vanishing Order.}   Let $X$ be a compact complex algebraic manifold of complex
dimension $n$ which is of general type.  Let
$$
\Phi=\sum_{m=1}^\infty\varepsilon_m\sum_{j=1}^{q_m}\left|s^{(m)}_j\right|^{\frac{2}{m}},
$$
where
$$
s^{(m)}_1,\cdots,s^{(m)}_{q_m}\in\Gamma\left(X,\,mK_X\right)
$$
form a basis over ${\mathbb C}$ and $\varepsilon_m>0$ approaches $0$
so fast as $m\to\infty$ that locally the infinite series which
defines $\Phi$ converges uniformly.   The metric 
defined by $\frac{1}{\Phi}$ for $K_X$ is called {\it a metric of minimum singularity}, which we also sometimes denote by $e^{-\varphi}$.

\medbreak For $N\in{\mathbb N}$ the function
$$
\Phi_N=\sum_{m=1}^N\varepsilon_m\sum_{j=1}^{q_m}\left|s^{(m)}_j\right|^{\frac{2}{m}}
$$
is called the $N$-th {\it truncation} of $\Phi$.   For a hypersurface $Y$ in $X$ and a regular point $P_0$ of $Y$ the {\it stable vanishing order} across $Y$ is at $P_0$ the infiumum of the vanishing order of the local multi-valued holomorphic function $\left(s^{(m)}_j\right)^{\frac{1}{m}}$ across $Y$ at $P_0$ for all $m\in{\mathbb N}$ and $1\leq j\leq q_m$.   The {\it generic stable vanishing order} across $Y$ is the stable vanishing order across $Y$ at a generic point $P_0$ of $Y$.

\medbreak Let $V$ be a subvariety of $X$ and $P_0$ be a regular point of $V$. We say that the {\it stable vanishing order is achieved} across $V$ at $P_0$ if there exist
some positive number $C_V$ and some $m_{P_0}\in{\mathbb N}$ such that $\Phi\leq C_V\Phi_{m_{P_0}}$  on $U_{P_0}$ for some open neighborhood $U_{P_0}$ of $P_0$ in $X$.  We say that the {\it generic stable vanishing order is achieved} across $V$ if the stable vanishing order is achieved across $V$ at some point $P_0$ of $V$.
We say that {\it all the stable vanishing orders are precisely achieved} if there exists some $m_0\in{\mathbb N}$ such that
$\Phi\leq C\Phi_{m_0}$
on $X$ for some positive constant $C$.

\medbreak\noindent(2.5.1) {\it Rationality of Stable Vanishing Orders.}   In many steps of the analytic proof of the finite generation of the canonical ring the rationality of every stable vanishing order needs to be verified so that the ${\mathbb R}$-divisor $\alpha Y$ can be regarded as a ${\mathbb Q}$-line bundle when $\alpha$ is the generic stable vanishing order across a hypersurface $Y$.  Such a verification of the rationality of every generic stable vanishing order $\alpha$ comes from the fact that the arguments in the analytic proof of the finite generation of the canonical ring will lead to a multi-valued holomorphic section giving a smaller generic vanishing order across $Y$ than $\alpha$ if $\alpha$ is irrational.  In this article we will not put in the verification for the rationality of the stable vanishing orders and will assume that all the stable vanishing orders occurring in the arguments presented in this article to be rational.

\bigbreak\noindent(2.6) {\it Finite Generation of Canonical Ring
From Achievement of Stable Vanishing Order.} Suppose all the stable
vanishing orders are achieved at every point of $X$ by the $m_0$-th truncation $\Phi_{m_0}$ of $\Phi$ for some $m_0\in{\mathbb N}$ so that $\Phi\leq C\Phi_{m_0}$ for some constant $C$. Denote $\left(m_0\right)!$ by $m_1$. Then the
canonical ring
$$
\bigoplus_{m=1}^\infty\Gamma\left(X, mK_X\right)
$$
is generated by
$$
\bigoplus_{m=1}^{\left(n+2\right)m_1}\Gamma\left(X, mK_X\right)
$$
and hence is finitely generated by the finite set of elements
$$
\left\{s^{(m)}_j\right\}_{1\leq m\leq(n+2)m_1,\,1\leq j\leq q_m}.
$$

\medbreak\noindent{\it Proof.}  Since $m_1=\left(m_0\right)!$ is divisible by $\nu$ and $\left(s^{(\nu)}_j\right)^{\frac{m_1}{\nu}}\in\Gamma\left(X, m_1K_X\right)$ for $1\leq\nu\leq m_0$, it follows that
$$\left(\sum_{\nu=1}^{m_0}{\sum_{j=1}^{q_\nu}}\left|s^{(\nu)}_j\right|^{\frac{2}{\nu}}\right)^{m_1}\leq C_1\sum_{j=1}^{q_{m_1}}\left|s^{\left(m_1\right)}_j\right|^2$$ for some constant $C_1$ and
$$\left(\Phi_{m_0}\right)^{m_1}\leq C_2\sum_{j=1}^{q_{m_1}}\left|s^{\left(m_1\right)}_j\right|^2$$ for some constant $C_2$.  
Let $e^{-\varphi}=\frac{1}{\Phi}$.
For $m>(n+2)m_1$ and any $s\in\Gamma\left(X, mK_X\right)$ we have
$$
\int_X\frac{\left|s\right|^2e^{-\left(m-\left(n+2\right)m_1-1\right)\varphi}}
{\left(\sum_{j=1}^{q_{m_1}}\left|s^{\left(m_1\right)}_j\right|^2\right)^{n+2}}<\infty,
$$
because $\left|s\right|^2\leq C_3\Phi^m$ on $X$ for some constant $C_3$ and $\Phi\leq C\Phi_{m_0}$. By Skoda's theorem on ideal generation ((2.3) and (2.4)) there
exist
$$
h_1,\cdots,h_{q_{m_1}}\in\Gamma\left(X,\left(m-m_1\right)K_X\right).  
$$
such that $s=\sum_{j=1}^{q_{m_1}}h_j s^{\left(m_1\right)}_j$.   If
$m-\left(n+2\right)m_1$ is still greater than $\left(n+2\right)m_1$,
we can apply the argument to each $h_j$ instead of $s$ until we get
$$
s\in\Gamma\left(X, \left(m-\ell m_1\right)K_X\right)\left(\Gamma\left(X, m_1 K_X\right)^\ell\right),
$$
where $\ell$ is the smallest integer such that $m-\ell m_1<(n+2)m_1$.  Thus $s$ is generated by
$$
\left\{s^{(m)}_j\right\}_{1\leq m\leq(n+2)m_1,\,1\leq j\leq q_m}.
$$
Q.E.D.

\bigbreak\noindent{\sc \S3. Decomposition of Closed Positive (1,1)-Currents and their Modified Restrictions to Hypersurfaces .}

\bigbreak\noindent(3.1) {\it Lelong Numbers of Closed 
Positive (1,1)-Current.}  For a closed positive $(1,1)$-current $\Theta$ on some
open subset $G$ of ${\mathbb C}^n$, the {\it Lelong number} of
$\Theta$ at a point $P_0$ of $G$ is the limit of
$$
\frac{\ \ \int_{B_n(P_0,r)}{\rm trace\,}\Theta\ \ }{{\rm
Vol\,}\left(B_{n-1}(0,r)\right)}
$$
as $r\to 0$, where $B_m\left(Q,r\right)$ is the open ball in
${\mathbb C}^m$ of radius $r$ centered at $Q$ and ${\rm
Vol\,}\left(B_m\left(Q,r\right)\right)$ is its volume and ${\rm
trace\,}\Theta$ is
$$\Theta\wedge\frac{1}{(n-1)!}\left(\sum_{j=1}^n\frac{\sqrt{-1}}{2}\,dz_j\wedge
d\overline{z_j}\right)^{n-1}$$ with $z_1,\cdots,z_n$ being the
coordinates of ${\mathbb C}^n$.   For $c>0$ the set $E_c=E_c\left(\Theta\right)$ consisting of all points of $G$ where the Lelong number of $\Theta$ is $\geq c$ is a complex-analytic subset of $G$ and is called a {\it Lelong set} of $\Theta$ (see, for example, [Lelong 1968, Siu 1974]).

\medbreak For a complex hypersurface $Y$ in $G$, integration over the regular points of $Y$ defines a closed positive $(1,1)$-current on $G$, which we denote by $\left[Y\right]$ or simply by $Y$.
The Lelong number of
$\left[Y\right]$ at a point $P_0$ of $Y$ is the multiplicity of $Y$ at $P_0$.

\bigbreak\noindent(3.2) {\it Canonical Decomposition of Closed
Positive (1,1)-Current.} Let $\Theta$ be a closed positive
$(1,1)$-current on a complex manifold $X$.  Then $\Theta$ admits a
unique {\it decomposition} of the following form
$$ \Theta=\sum_{j=1}^J\gamma_j\left[V_j\right]+R,
$$
where $\gamma_j>0$, $J\in{\mathbb N}\cup\left\{0,\infty\right\}$,
$V_j$ is a complex hypersurface in $X$ and the Lelong number of the
remainder $R$ is zero outside a countable union of subvarieties of
codimension $\geq 2$ in $X$ [Siu 1974].  We consider the {\it dichotomy}
into two cases.  The first case is either $R\not=0$ or $J=\infty$.
The second case is when both $R=0$ and $J$ is finite.

\bigbreak\noindent(3.3) {\it Modifications of Restrictions of Closed Positive (1.1)-Currents to Hypersurfaces.}  Let $X$ be a compact complex algebraic
manifold of general type and let $e^{-\varphi}=\frac{1}{\Phi}$ be
the metric of minimum singularity as defined in (2.5).  Let
$$
\Theta_\varphi=\frac{\sqrt{-1}}{2\pi}\partial\bar\partial\varphi
$$
be the curvature current of the metric $e^{-\varphi}$ of $K_X$.  Let
$Y$ be a nonsingular hypersurface in $X$.   Let $\gamma$ be the generic stable
vanishing order across $Y$, which is equal to Lelong number of
$\Theta_\varphi$ at a generic point of $Y$.

\medbreak We are going to define the restriction
to $Y$ of the closed positive $(1,1)$-current
$\Theta-\gamma\left[Y\right]$.  We call such a restriction to $Y$
the {\it modified restriction} of $\Theta$ to $Y$, because we are
restricting $\Theta$ after we modify it by subtracting
$\gamma\left[Y\right]$ from it.  
For $k\in{\mathbb N}$ let
$\gamma_k$ be the infimum of the vanishing order of the multi-valued
section $\left(s^{(m)}_j\right)^{\frac{1}{m}}$ across $Y$ for $1\leq
m\leq k$ and $1\leq j\leq q_m$. Consider the metric
$$
\frac{1}{\sum_{m=1}^k\varepsilon_m\sum_{j=1}^{q_m}
\left|\frac{\left(s^{(m)}_j\right)^{\frac{1}{m}}}{s_Y^{\gamma_k}}\right|^2}
$$
of the ${\mathbb Q}$-line bundle $\left.\left(L-\gamma_k
Y\right)\right|_Y$ on $Y$ and its curvature current
$$
\Theta_k=\frac{\sqrt{-1}}{2\pi}\partial\bar\partial\log\sum_{m=1}^k\varepsilon_m\sum_{j=1}^{q_m}
\left|\frac{\left(s^{(m)}_j\right)^{\frac{1}{m}}}{s_Y^{\gamma_k}}\right|^2
$$ which is a
closed positive $(1,1)$-current on $Y$.  We know that the sequence
$\gamma_k$ is non-increasing and its limit is $\gamma$ as
$k\to\infty$ so that the ${\mathbb Q}$-line bundle
$\left.\left(L-\gamma_k Y\right)\right|_Y$ on $Y$ approaches the
${\mathbb Q}$-line bundle $\left.\left(L-\gamma Y\right)\right|_Y$
on $Y$ as $k\to\infty$.  The restriction of the closed positive
$(1,1)$-current $\Theta-\gamma\left[Y\right]$ to $Y$ can be defined
as the (weak) limit of $\Theta_k$ (or its subsequence) as
$k\to\infty$.

\medbreak We do not consider using the metric
$$
\frac{1}{\sum_{m=1}^\infty\varepsilon_m\sum_{j=1}^{q_m}
\left|\frac{\left(s^{(m)}_j\right)^{\frac{1}{m}}}{s_Y^\gamma}\right|^2}
$$
of the ${\mathbb Q}$-line bundle $\left(L-\gamma Y\right)|_Y$ and then
using its curvature current as the restriction of the closed positive
$(1,1)$-current $\Theta-\gamma\left[Y\right]$ to $Y$, because
the vanishing order of $s^{(m)}_j$ across $Y$ may be strictly higher than 
$m\gamma$ for each $m\in{\mathbb N}$ and each $1\leq
j\leq q_m$ so that the multi-valued fraction
$$
\frac{\left(s^{(m)}_j\right)^{\frac{1}{m}}}{s_Y^\gamma}
$$
may still be identically zero on $Y$ for each $m\in{\mathbb N}$ and
each $1\leq j\leq q_m$ and the sum
$$
\sum_{m=1}^\infty\varepsilon_m\sum_{j=1}^{q_m}
\left|\frac{\left(s^{(m)}_j\right)^{\frac{1}{m}}}{s_Y^\gamma}\right|^2
$$
may be identically zero on $Y$, making it impossible to use such a definition. 

\medbreak The nonsingular hypersurface $Y$ in $X$ is said to belong to the first case of the dichotomy (respectively the second case of the dichotomy) when the closed positive $(1,1)$-current $\Theta-\gamma\left[Y\right]$ on $Y$ belongs to the first case of the dichotomy (respectively the second case of the dichotomy).

\medbreak\noindent(3.3.1) {\it Remark.}  When $\Theta-\gamma\left[Y\right]$ is in the second case of the dichotomy in the sense of (3.2), if $\sigma$ is a nonzero element of $\Gamma\left(X, mK_X\right)$ whose generic vanishing order across $Y$ is $m\gamma$, then for any other nonzero element $\hat\sigma$  of $\Gamma\left(X, \hat mK_X\right)$ whose generic vanishing order across $Y$ is $\hat m\,\gamma$ the quotient $\frac{\sigma^{\hat m}}{\hat\sigma^m}$ has to be equal to a nonzero constant on $Y$.  However, in general it does not mean that the line bundle $K_X-\gamma Y$ over $Y$ is flat.  The limitation is only on elements of $\Gamma\left(Y, m\left(K_X-\gamma Y\right)\right)$ which are extendible to elements of $\Gamma\left(X, m\left(K_X-\gamma Y\right)\right)$.

\bigbreak\noindent{\sc \S4.  Discrepancy Subspaces} 

\bigbreak\noindent(4.1) {\it Motivation for the Definition of Discrepancy Subspaces.}  As explained in  Item (H) in the Introduction, we seek to describe the deviation of $mK_X$ from a prescribed sufficiently ample line bundle $A$ so that the deviation is formulated to be independent of $m$.  This is essential to terminate the infinite number of blow-ups in the process of proving that all the stable vanishing orders are precisely achieved.   The final goal is, of course, to show that $\Phi_m$ is comparable to $\Phi$ for some $m\in{\mathbb N}$.   We would like to measure the failure of the comparability between $\Phi_m$ and $\Phi$ in a way which is independent of $m$.  The deviation, uniform in $m$, for $mK_X$ from $A$ is just an intermediate step to achieve the final goal of measuring the failure of the comparability between $\Phi_m$ and $\Phi$ in a way which is independent of $m$.  A discrepancy subspace is to precisely formulate this measurement of uniform failure of comparability.

\medbreak In order to know why the definition which we are going to describe is to be formulated in the way given here, we have to go back to \S2 where it is explained how the analytic proof of the finite generation of the canonical ring is reduced to the precise achievement of all stable vanishing orders.  For the purpose of invoking Skoda's theorem on ideal generation to show that $s^{(m)}_j$ is expressible in terms of $s^{\left(m_0\right)}_1,\cdots,s^{\left(m_0\right)}_{q_m}$ for $m\geq(n+2)m_0$ with coefficients in $\Gamma\left(X,\left(m-m_0\right)K_X\right)$, we would like to have
$$
\int_X\frac{\left|s^{\left(m\right)}_j\right|^2}
{\left(\sum_{j=1}^{q_{m_0}}\left|s^{\left(m_0\right)}_j\right|^2\right)^{n+2}\Phi^{m-m_0\left(n+2\right)}}<\infty.\leqno{(4.1.1)}
$$
The original philosophy of the use of multiplier ideal sheaves, both from Kohn's perspective of partial differential equations [Kohn 1979] and from Nadel's perspective of instability [Nadel 1990], is to introduce multipliers into the most crucial estimates when such estimates fail to hold (see the Appendix at the end of this article).   Here the crucial estimate is (4.1.1).  In case it fails to hold, we introduce a multiplier $f$, which is a holomorphic function germ, such that
$$\left|f\right|^2\frac{\left|s^{\left(m\right)}_j\right|^2}
{\left(\sum_{j=1}^{q_{m_0}}\left|s^{\left(m_0\right)}_j\right|^2\right)^{n+2}\Phi^{m-m_0\left(n+2\right)}}\leqno{(4.1.2)}
$$
is locally integrable for all $m\geq(n+2)m_0$ or for all $m\geq\hat m$ for some $\hat m$ independent of $m$.  Since
$$
\left|s^{\left(m\right)}_j\right|^2\leq\frac{1}{\left(\varepsilon_m\right)^m}\Phi^m,
$$
for the integrability of (4.1.2) it suffices to assume that
$$\left|f\right|^2\Phi^{m_0\left(n+2\right)}\leq C_{m_0}\left(\sum_{j=1}^{q_{m_0}}\left|s^{\left(m_0\right)}_j\right|^2\right)^{n+2}\leqno{(4.1.3)_{m_0}}
$$ for some constant $C_{m_0}$.   Here we switch from integrability in terms of the $L^2$ norm to the supremum norm, because the supremum norm is easier to keep track of when we blow up.  For the $L^2$ norm every time we blow up we have to worry about the contribution from the Jacobian determinants.  It simply makes things a bit more complicated.  See Remark (4.2.3).

\medbreak We do not know beforehand which $m_0$ is able to give finite generation.  Our goal is to increase the ideal formed by the multipliers $f$ until it becomes the unit ideal and the method is to increase the $m_0$ used.  Thus we need to have a multiplier $f$ which works for all $m_0$ sufficiently large.  In other words, we should consider multipliers $f$ such that $(4.1.3)_{m_0}$ holds for all $m_0$ sufficiently large.   The power of $n+2$ on both sides of $(4.1.3)_{m_0}$ can be removed at the expense of using the normalization of raising the multiplier ideal sheaves to the $(n+2)$-th power at the end.   Note that this kind of multiplier ideal sheaves is {\it dynamic}, because a sequence of inequalities is used in the definition instead of just one single inequality.  As remarked in the Introduction this is more in keeping with the original notions of multiplier ideal sheaves of Kohn and Nadel and is different from the {\it static} kind which is usually used in algebraic geometry and which is defined by a single inequality.

\medbreak The above discussion motivates us to introduce a coherent ideal sheaf ${\mathcal J}$ on $X$ such that the inequality
$$
\left|s_{\mathcal J}\right|^2\,\Phi^m\leq\tilde C_{m,{\mathcal
J}}\sum_{j=1}^{q_m}\left|s^{(m)}_j\right|^2\leqno{(4.1.4)}
$$
holds on $X$ for some constant $\tilde C_{m,{\mathcal J}}$, where the notation $\left|s_{\mathcal J}\right|^2$ means the
following. For a coherent ideal sheaf ${\mathcal I}$ on $X$
generated locally by holomorphic function germs
$\tau_1,\cdots,\tau_\ell$, we define
$$
\left|s_{\mathcal I}\right|^2=\sum_{j=1}^\ell\left|\tau_j\right|^2.
$$
Now the positive constants $\varepsilon_\ell$ used in the definition of
$$
\Phi=\sum_{\ell=1}^\infty\varepsilon_\ell\sum_{j=1}^{q_\ell}\left|s^{(
\ell)}_j\right|^{\frac{2}{\ell}}
$$
are quite arbitrarily chosen as long as the sequence which they form decrease sufficiently rapidly to guarantee the convergence of the infinite series.  In the inequality (4.1.4) we could allow the set of positive constants $\left\{\varepsilon_\ell\right\}$ for $\Phi$ on the left-hand side to depend on $m$ so that $\left\{\varepsilon_\ell\right\}$ is replaced by $\left\{\varepsilon_{m,\ell}\right\}$.  In
other words, instead of (4.1.4) we use the inequality
$$
\left|s_{\mathcal J}\right|^2\left(\check{\Phi}_m\right)^m\leq\tilde
C_{m,{\mathcal J}}\sum_{j=1}^{q_m}\left|s^{(m)}_j\right|^2,\leqno{(4.1.5)}
$$
where
$$
\check{\Phi}_m=\sum_{\ell=1}^\infty\varepsilon_{\ell,m}\sum_{j=1}^{q_\ell}\left|s^{(\ell)}_j\right|^{\frac{2}{\ell}}
$$
for some positive constants $\varepsilon_{\ell,m}$.  We now incorporate the constants $\varepsilon_{\ell,m}$ and $\tilde
C_{m,{\mathcal J}}$ together to give the following definition of a discrepancy subspace.

\bigbreak\noindent(4.2) {\it Definition of Discrepancy Subspace.}
Let ${\mathcal J}$ be a coherent ideal sheaf on $X$. The stable
vanishing order of the canonical line bundle of $X$ is said to be
{\it precisely achieved modulo the subspace of $X$ defined by
${\mathcal J}$} if there exist some positive integer $m_{\mathcal
J}$ and some positive constant $C_{m,k,{\mathcal J}}$ for $k,
m\in{\mathbb N}$ with $m\geq m_{\mathcal J}$ such that the
inequality
$$
\left|s_{\mathcal
J}\right|^2\,\sum_{j=1}^{q_k}\left|s^{(k)}_j\right|^{\frac{2m}{k}}\leq
C_{m,k,{\mathcal
J}}\sum_{j=1}^{q_m}\left|s^{(m)}_j\right|^2\leqno{(4.2.1)}
$$
holds on $X$ for all $k, m\in{\mathbb N}$ with $m\geq m_{\mathcal
J}$.  Let $Z$ be the zero-set of the coherent ideal sheaf ${\mathcal J}$.
We call the ringed space $\left(Z,{\mathcal O}_X\left/{\mathcal
J}\right.\right)$ a {\it discrepancy subspace}.  We call the
coherent ideal sheaf ${\mathcal J}$ a {\it discrepancy ideal sheaf}.

\medbreak\noindent(4.2.2){\it Remark.} An intuitive way of describing the discrepancy subspace is that it gives an $m$-independent bound for the difference between any stable vanishing order defined by {\it multi-valued} global sections of the $m$-canonical line bundle and the vanishing order actually achieved by global {\it single-valued} $m$-canonical sections, as described by the sum of absolute-value-squares of local generators of the ideal sheaves. 

\medbreak\noindent(4.2.3) {\it Remark.}  In a blow-up $\tilde X\to X$ of $X$, the
r\^ole played by the adjunction formula is canceled by its effect on
both sides of $(4.2.1)$ and $\left|s_{\mathcal J}\right|^2$ simply
transforms as a lifting of a local function from $X$ to $\tilde X$.
This enables us to assume that, after replacing $X$ by its blowup,
${\mathcal J}$ is the ideal sheaf of a divisor whose components are
in normal crossing.  With the blow-up, we can use the technique of
the minimum center of log canonical singularity [Kawamata 1985,
Shokurov 1985]  or its analogue in our case.

\medbreak\noindent(4.2.4) {\it Remark.}  A consequence of the inequality (4.2.1) which defines the discrepancy subspace is that 
$$
{\mathcal J}{\mathcal I}_{m\log\Phi}\subset{\mathcal
I}_{\log\sum_{j=1}^m\left|s^{(m)}_j\right|^2},
$$
which means that the {\it conductor} ${\mathcal J}$ from the ideal sheaf ${\mathcal I}_{m\log\Phi}$ to the ideal sheaf ${\mathcal I}_{\log\sum_{j=1}^m\left|s^{(m)}_j\right|^2}$ is independent of $m$.

\medbreak\noindent(4.2.5) {\it Remark.}  From (4.2.1) it follows that after replacing $X$ an appropriate blow-up whose base-point set is a union of nonsingular hypersurfaces in normal crossing, the discrepancy ideal sheaf ${\mathcal J}$ can be replaced by the multiplier ideal sheaf of the metric
$$
e^{-\psi_{\mathcal J}}=\frac{\prod_{\nu=1}^\ell\left|s_{Y_\nu}\right|^{2p\alpha_\nu}}{\sum_{j=1}^k\left|\sigma_j\right|^2}
$$
of $p\left(K_X-\sum_{\nu=1}^\ell\alpha_\nu Y_\nu\right)$ for some sufficiently large $p\in{\mathbb N}$, where $\left\{Y_\nu\right\}_{\nu=1}^k$ is the set of all hypersurfaces of the base-point set not contained in the zero-set of ${\mathcal J}$ and $\alpha_\nu$ is the generic stable vanishing order across $Y_\nu$ and $\sigma_1,\cdots,\sigma_k$ are some multi-valued holomorphic sections of $pK_X$.

\bigbreak\noindent(4.3) {\it Intersection of Discrepancy Subspaces.}  The definition of discrepancy subspaces given in (4.2) allows us
to take the intersection of two discrepancy subspaces in order to decrease a discrepancy subspace all the way down to the empty set.  Suppose we have two discrepancy subspaces defined by discrepancy ideal sheaves ${\mathcal
J}$ and $\tilde{\mathcal
J}$ so that
$$
\left|s_{\mathcal
J}\right|^2\,\sum_{j=1}^{q_k}\left|s^{(k)}_j\right|^{\frac{2m}{k}}\leq
C_{m,k,{\mathcal
J}}\sum_{j=1}^{q_m}\left|s^{(m)}_j\right|^2
$$
on $X$ for all $k, m\in{\mathbb N}$ with $m\geq m_{\mathcal J}$ and 
$$
\left|s_{\tilde{\mathcal
J}}\right|^2\,\sum_{j=1}^{q_k}\left|s^{(k)}_j\right|^{\frac{2m}{k}}\leq
C_{m,k,\tilde{\mathcal
J}}\sum_{j=1}^{q_m}\left|s^{(m)}_j\right|^2$$
on $X$ for all $k, m\in{\mathbb N}$ with $m\geq m_{\tilde{\mathcal
J}}$.  Then we can define their intersection discrepancy subspace with discrepancy ideal sheaf ${\mathcal K}$  by setting ${\mathcal K}$ as the sum of ${\mathcal
J}$ and $\tilde{\mathcal J}$, $m_{\mathcal K}$ as the maximum
of $m_{\mathcal J}$ and $m_{\tilde{\mathcal J}}$, and 
$C_{m,k,{\mathcal K}}$ as the sum of $C_{m,k,{\mathcal J}}$ and
$C_{m,k,{\tilde{\mathcal J}}}$ so that
$$
\left|s_{\mathcal
K}\right|^2\,\sum_{j=1}^{q_k}\left|s^{(k)}_j\right|^{\frac{2m}{k}}\leq
C_{m,k,{\mathcal
K}}\sum_{j=1}^{q_m}\left|s^{(m)}_j\right|^2
$$
on $X$ for all $k, m\in{\mathbb N}$ with $m\geq m_{\mathcal K}$.  

\bigbreak\noindent(4.3.1) {\it Remark.}  If we do not use the definition of discrepancy subspaces given in (4.2) and choose, instead, to define a discrepancy subspace ${\mathcal J}_{m_0}$ for a fixed $m_0$-canonical line bundle by
$$
\left|s_{{\mathcal
J}_{m_0}}\right|^2\,\sum_{j=1}^{q_k}\left|s^{(k)}_j\right|^{\frac{2m_0}{k}}\leq
C_{k,{\mathcal
J}_{m_0}}\sum_{j=1}^{q_{m_0}}\left|s^{\left(m_0\right)}_j\right|^2\leqno{(4.3.1.1)}
$$
on $X$ for all $k\in{\mathbb N}$ with $k\geq\kappa_{m_0,{\mathcal J}_{m_0}}$, motivated by
$(4.1.3)_{m_0}$, then when for $\tilde m_0=pm_0$ with an integer $p>1$ there is another discrepancy subspace $\tilde{\mathcal J}_{\tilde m_0}$ for a fixed $\tilde m_0$-canonical line bundle defined by
$$
\left|s_{\tilde{\mathcal
J}_{\tilde m_0}}\right|^2\,\sum_{j=1}^{q_k}\left|s^{(k)}_j\right|^{\frac{2\tilde m_0}{k}}\leq
\tilde C_{k, \tilde{\mathcal
J}_{\tilde m_0}}\sum_{j=1}^{q_{\tilde m_0}}\left|s^{\left(\tilde m_0\right)}_j\right|^2
$$
on $X$ for all $k\in{\mathbb N}$ with $k\geq\tilde\kappa_{\tilde m_0,\tilde {\mathcal J}_{\tilde m_0}}$, we can only get the inequality
$$
\left(\left|s_{{\mathcal
J}_{m_0}}\right|^{2p}+\left|s_{\tilde{\mathcal
J}_{\tilde m_0}}\right|^2\right)\,\sum_{j=1}^{q_k}\left|s^{(k)}_j\right|^{\frac{2m}{k}}\leq
\hat C_{k, \tilde{\mathcal
J}_{\tilde m_0}}\sum_{j=1}^{q_{\tilde m_0}}\left|s^{\left(\tilde m_0\right)}_j\right|^2
$$
on $X$ for all $k\in{\mathbb N}$ with $k\geq\max\left(\kappa_{m_0,{\mathcal J}_{m_0}},\tilde\kappa_{\tilde m_0,\tilde {\mathcal J}_{\tilde m_0}}\right)$.  With this kind of definition for only a fixed pluricanonical line bundle, we can only replace ${\mathcal J}_{m_0}$ by the smaller ideal sheaf $\left({\mathcal J}_{m_0}\right)^p+\tilde{\mathcal J}_{\tilde m_0}$ instead of by the ideal sheaf ${\mathcal J}_{m_0}+\tilde{\mathcal J}_{\tilde m_0}$ as in (4.3).  Only by replacing ${\mathcal J}_{m_0}$ by ${\mathcal J}_{m_0}+\tilde{\mathcal J}_{\tilde m_0}$ can we decrease the discrepancy subspace to the empty set.

\medbreak In contrast to the definition given in (4.2) this definition for only a fixed pluricanonical line bundle given by (4.3.1.1) involves only one inequality when $\Phi^{m_0}$ with another sequence of rapidly decreasing positive coefficients $\left\{\varepsilon_m\right\}_{m\in{\mathbb N}}$ and is {\it static}, whereas the definition given in (4.2) even when it is expressed in terms of powers of $\Phi$ involves a sequence of inequalities and is {\it dynamic}.

\medbreak To prepare for the construction of discrepancy subspaces we need the following lemma from the round-up and round-down properties of the usual multiplier ideal sheaves in algebraic geometry.

\bigbreak\noindent(4.4) {\it Lemma on Sup Norm Domination of Metric
by Generators of Multiplier Ideal.} Let $f_j$ be holomorphic
functions on some open neighborhood $U$ of the origin in ${\mathbb
C}^n$. Let $\varepsilon_j>0$ and $m_j\in{\mathbb N}$ so that
$$
\Psi=\sum_{j=1}^\infty\varepsilon_j\left|f_j\right|^{\frac{2}{m_j}}
$$
converges uniformly on compact subsets of $U$.  Let ${\mathcal J}$
be the multiplier ideal sheaf of the metric $\frac{1}{\Psi}$ and
$g_1,\cdots,g_\ell$ be holomorphic function germs on ${\mathbb C}^n$
at the origin such that the stalk of ${\mathcal J}$ at the origin is
generated by $g_1,\cdots,g_\ell$ over ${\mathcal O}_{{\mathbb
C}^n,0}$. Then there exists an open neighborhood $W$ of the origin
in ${\mathbb C}^n$ where $g_1,\cdots,g_\ell$ are defined and there
exists a positive constant $C_j$ such that
$$
\left|f_j\right|^{\frac{2}{m_j}}\leq
C_j\sum_{k=1}^\ell\left|g_k\right|^2
$$
on $W$.

\medbreak\noindent(4.4.1) {\it Remark.}  The geometric reason for
this lemma is that the minimum of the orders of the zeros of the
generators of a multiplier ideal ${\mathcal J}$ should be no more
than the order of the pole of the metric $\frac{1}{\Psi}$.  A proof,
for example, is given in [Proposition 3.1, Demailly 1992].

\bigbreak\noindent(4.5) {\it Construction of Initial Codimension-One
Discrepancy Subspace.}  As the first step we now construct the
initial codimension-one discrepancy subspace by using
the technique of the global generation of the multiplier ideal sheaf
(2.2) and the decomposition of $K_X$ as a sum of an ample
${\mathbb Q}$-line bundle and an effective ${\mathbb Q}$-divisor
from the general type property of $X$.

\medbreak Let $A$ be an ample line bundle on $X$ which is ample
enough for the global generation of multiplier ideal sheaves, as described in (2.2).
We write $aK_X=A+D$, where $D$ is an effective
divisor in $X$ and $a$ is a positive integer. We use the metric
$$\frac{1}{\Phi^m\left|s_D\right|^2}$$ for the line bundle
$$mK_X+D=\left(m+a\right)K_X-A.$$   Let ${\mathcal I}^{(m)}$ be the multiplier ideal
sheaf of the metric
$$\frac{1}{\Phi^m\left|s_D\right|^2}.$$   Then the multiplier ideal sheaf
${\mathcal I}^{(m)}$ is generated by elements of
$$\displaylines{\qquad\qquad\qquad\Gamma\left(X,{\mathcal I}^{(m)}\left(mK_X+D+A+K_X\right)\right)\hfill\cr\hfill=
\Gamma\left(X,{\mathcal
I}^{(m)}\left(\left(m+a+1\right)K_X\right)\right)\subset
\Gamma\left(X,\left(m+a+1\right)K_X\right).\qquad\cr}$$  From 
Lemma (4.4) on the sup norm domination of a metric by the generators of its multiplier
ideal we conclude that
$$
\left|s^{(k)}_j\right|^{\frac{2m}{k}}\left|s_D\right|^2\leq
C_{k,j,m}\sum_{j=1}^{q_{m+a+1}}\left|s^{(m+a+1)}_j\right|^2
$$
for $k\in{\mathbb N}$, which implies
$$
\left|s^{(k)}_j\right|^{\frac{2(m+a+1)}{k}}\left|s_D\right|^2\leq
\tilde C_{k,j,m}\sum_{j=1}^{q_{m+a+1}}\left|s^{(m+a+1)}_j\right|^2
$$
for $k\in{\mathbb N}$, because each $s^{(k)}_j$ is a local holomorphic function, where $C_{k,j,m}$ and $\tilde C_{k,j,m}$ are constants.
This shows that we can choose ${\mathcal J}$
to be the ideal sheaf generated by $s_D$ and choose $m_{\mathcal J}$
as $a+2$.

\medbreak\noindent(4.5.1) {\it Remark on Key Points in Construction of Discrepancy Subspace.}  The three key points of the construction of the initial discrepancy subspace of codimension one in (4.5) are the following.
\begin{itemize}\item[(i)]  Kodaira's trick of squeezing out some ample-line-bundle part of $K_X$ by writing $K_X=A+D$ for some effective ${\mathbb Q}$-divisor $D$ and some ample ${\mathbb Q}$-line bundle $A$, because of the growth order of the dimension of the space of all global holomorphic sections of $mK_X$ as a function of $m$ so that we can find some nontrivial global $m$-canonical section vanishing on the divisor of an ample line bundle.  This is a matter of the growth order of the dimension of the space of global holomorphic $m$-canonical sections.
\item[(ii)]  The use of a high enough multiple of the sqeezed-out ample line bundle $A$ to globally generate any multiplier ideal sheaf (with the multiple independent of the multiplier ideal sheaf). 
\item[(iii)] The relation between the local supremum norm and the local $L^2$ norm used in defining multiplier ideal sheaves converts the global generation of multiplier ideal sheaves with a fixed sufficiently ample twisting into the inequalities in the definition of discrepancy subspaces.
\end{itemize}
This argument of the above three points can be applied also to the case of a fibration when the holomorphic sections on the fibration are constant along generic fibers and the ample line bundle is only on the base of the fibration and the restrictions of the multiplier ideal sheaves to generic fibers are already known to be globally generated without any ample twisting.

\medbreak Such a fibration is needed, because in order to decrease the discrepancy subspace we have to replace $X$ by the discrepancy subspace in our argument and the discrepancy subspace in general is not of general type, but can be described by a fibration whose generic fibers belong to the second case of the dichotomy.

\bigbreak\noindent(4.6) {\it Decreasing Discrepancy Subspace by
Applying the Argument of Constructing Initial Codimension-One
Discrepancy Subspace to Fibrations.}   In order to decrease the initial codimension-one discrepancy subspace constructed in (4.5) we have to apply the argument of (4.5) to a fibration on the initial codimension-one discrepancy subspace $\left(D,{\mathcal O}_X\left/{\mathcal J}\right)\right.$ as briefly indicated at the end of (4.5.1).  Since the techniques require the extension of sections for the two cases of the dichotomy, before we can decrease the discrepancy subspace we have to first explain the extension techniques of sections for the two cases of the dichotomy.

\bigbreak\noindent{\sc \S5 Construction of Pluricanonical Sections with Fixed Sufficiently Ample Twisting}

\bigbreak\noindent(5.1)  For a compact complex algebraic manifold $X$ of general type, we consider a nonsingular hypersurface $Y$ in $X$ whose generic stable vanishing order is $\gamma$.   Our goal is to show that the generic stable vanishing order across $Y$ is actually achieved by some global $m$-canonical section $s\in\Gamma\left(X, mK_X\right)$ of $X$ in the sense that the vanishing order of $s$ across $Y$ at some generic point of $Y$ is precisely $m\alpha$.  

\medbreak In this \S5 we will not be able to produce right away such an $s\in\Gamma\left(X, mK_X\right)$.  We will assume that $Y$ belongs to the first case of the  dichotomy.  At any prescribed point $P_0$ of $X$ we will produce for an appropriately large $\hat p\in{\mathbb N}$ a strictly positively curved metric $\frac{1}{\Phi_{Y,P_0,\varepsilon, N}}$ of $\hat p K_X$ on $X$ defined by a finite sum $\Phi_{Y,P_0,\varepsilon, N}$ of absolute-value-squares of multi-valued holomorphic sections of $\hat p K_X$ over $X$ which has high prescribed vanishing order $N$ at  $P_0$ and yet has a vanishing order $\gamma$ across $Y$ at a generic point of $Y$ as close to the best expected value $2\hat p\alpha$ as any prescribed small error $<\varepsilon$ (see (5.5) below).   

\medbreak The prescribed small error comes from using the general type property of $X$ to fulfill the requirement of the strict positivity of the curvature current of the metric  $\frac{1}{\Phi_{Y,P_0,\varepsilon, N}}$ and to take care of some ample line bundle $A$ of $X$ used in an intermediate step to construct multi-valued holomorphic sections of $p\left(X-\alpha Y\right)+A$ over $Y$. 

\medbreak This metric $\frac{1}{\Phi_{Y,P_0,\varepsilon, N}}$ enables us to reduce the problem to producing a section on the subvariety $Y_1$ of $Y$ where the vanishing order of $\frac{\Phi_{Y,P_0,\varepsilon, N}}{\left|s_Y\right|^{2\gamma}}$ is high.  This kind of reduction is one of the techniques commonly used in the proof of non-vanishing theorems and in the study of problems related to the Fujita conjecture [Shokurov 1985, Kawamata 1985, Fujita 1987].  

\medbreak The reason why we cannot produce right away the element $s\in\Gamma\left(X, mK_X\right)$ just described is that eventually the process of replacing $Y$ by $Y_1$ when continued will lead to the second case of the dichotomy if $\alpha>0$, otherwise the definition of the generic stable vanishing order $\gamma$ will be contradicted (see (5.5.2) below for more details).  So the construction of the metric  $\frac{1}{\Phi_{Y,P_0,\varepsilon, N}}$ is the best we can do when we consider only the first case of the dichotomy.  This technique of the construction and the use of such a kind of metric constructed from the first case of dichotomy will also be applied to a generic fiber in a fibration which arises in decreasing the discrepancy subspace as mentioned in (4.6).  

\medbreak When eventually the second case of the dichotomy occurs, it will be handled in \S6 by using the section given by the decomposition of the modified restriction of the curvature current and the technique of Shokurov of using the theorem of Hirzebruch-Riemann-Roch to compare the arithmetic genus of the line bundle and
that of its twisting by a flat line bundle [Shokurov
1985].

\bigbreak\noindent(5.2) {\it Proposition (Sections of Amply Twisted Multiple of Line Bundle).}  Let $Y$ be a compact complex algebraic manifold of complex dimension $n$.   Let $A$ be a very ample holomorphic line bundle  on $Y$ with $A-K_Y$ also very ample.  Let $L$ be a holomorphic line bundle on $Y$ with metric $e^{-\varphi}$ whose curvature current $\Theta$ is a closed positive $(1,1)$-current in the first case of the dichotomy in the sense of (3.2), {\it i.e.,} in the decomposition
$$ \Theta=\sum_{j=1}^J\gamma_j\left[V_j\right]+R,
$$
with $\gamma_j>0$, $J\in{\mathbb N}\cup\left\{0,\infty\right\}$,
$V_j$ being a complex hypersurface and the Lelong number of the
remainder $R$ being zero outside a countable union of subvarieties of
codimension $\geq 2$ in $Y$, either  $R\not=0$ or $J=\infty$.  Then there exists some sequence $\left\{p_\nu\right\}_{\nu\in{\mathbb N}}$ of positive integers such that
$$
\lim_{\nu\to\infty}\dim_{\mathbb C}\Gamma\left(Y, {\mathcal I}_{p_\nu\varphi}\left(p_\nu L+n A\right)\right)=\infty.
$$

\medbreak\noindent{\it Proof.}  Let $p$ be a positive integer. 
Take $P_0\in Y$ and we will
also impose more conditions on $P_0$ later. Let $s_1$ be a generic
element of $\Gamma\left(Y,\,A\right)$ vanishing at $P_0$ so that the
sequence
$$
\displaylines{\qquad 0\to{\mathcal
I}_{p\varphi}\left(pL+A\right)\stackrel{\theta_{s_1}}{\longrightarrow}{\mathcal
I}_{p\varphi}\left(pL+2A\right)\hfill\cr\hfill\to\left({\mathcal
I}_{p\varphi}\left/s_1{\mathcal
I}_{p\varphi}\right.\right)\left(pL+2A\right)\to
0\qquad\cr}
$$
is exact, where $\theta_{s_1}$ is defined by multiplication by
$s_1$. For this step we have to make sure that the maximum ideal
${\mathfrak m}_{Y,P_0}$ of $Y$ at $P_0$ is not an associated prime
ideal in the primary decomposition of the stalk of the ideal sheaf
${\mathcal I}_{p\varphi}$ at $P_0$. This means that for each $p$ we
have to impose the condition that $P_0$ does not belong to some
finite subset $Z_0$ of $Y$.  Let $Y_1$ be the zero-set of $s_1$ and
$$
{\mathcal O}_{Y_1}=\left({\mathcal O}_Y\left/s_1{\mathcal
O}_Y\right.\right)|_{Y_1},
$$
which we can assume to be regular with ideal sheaf equal to
$s_1{\mathcal O}_Y$ because $s_1$ is generic element of
$\Gamma\left(Y,A\right)$ vanishing at $P_0$.  By choosing $s_1$
generically we can also assume that ${\mathcal
I}_{\left(p\varphi|_{Y_1}\right)}={\mathcal
I}_{p\varphi}\left/s_1{\mathcal I}_{p\varphi}\right.$. We use
$\chi\left(\cdot,\,\cdot\right)$ to denote the arithmetic genus
which means
$$
\chi\left(\cdot,\,\cdot\right)=\sum_{\nu=0}^\infty(-1)^\nu\dim_{\mathbb
C}H^\nu\left(\cdot,\,\cdot\right).
$$
From the long cohomology exact sequence of the above short exact
sequence we obtain
$$
\displaylines{\chi\left(Y,\,{\mathcal
I}_{p\varphi}\left(pL+2A\right)\right)=\cr
\chi\left(Y,\,{\mathcal
I}_{p\varphi}\left(pL+A\right)\right)+\chi\left(Y_1,\,{\mathcal
I}_{\left(p\varphi|_{Y_1}\right)}\left(pL+2A\right)|_{Y_1}\right).\cr}
$$
Since $A-K_Y$ is ample and $2A-K_{Y_1}=A-K_Y$ is also ample, by the theorem of Kawamata-Viehweg-Nadel
$$
\displaylines{H^\nu\left(Y,{\mathcal
I}_{p\varphi}\left(pL+kA\right)\right)=0\quad{\rm for\ }\nu\geq
1\ \ {\rm and}\ \ k=1,2,\cr H^\nu\left(Y_1,{\mathcal
I}_{\left(p\varphi|_{Y_1}\right)}\left(\left(pL+2A\right)|_{Y_1}\right)\right)=0\quad{\rm
for\ }\nu\geq 1.\cr}
$$
so that
$$
\displaylines{\Gamma\left(Y,\,{\mathcal
I}_{p\varphi}\left(pL+2A\right)\right)=\cr
\Gamma\left(Y,\,{\mathcal
I}_{p\varphi}\left(pL+A\right)\right)+\Gamma\left(Y_1,\,{\mathcal
I}_{\left(p\varphi|_{Y_1}\right)}\left(\left(pL+2\right)A\right)|_{Y_1}\right)\cr}
$$

\medbreak\noindent(5.2.1) {\it Slicing by Ample Divisors Down to a
Curve.} Instead of one single generic element $s\in\Gamma\left(Y,A\right)$,
we can choose generically
$$
s_1, \cdots, s_{n-1} \in\Gamma\left(Y,\,A\right)
$$
all vanishing at $P_0$ so that inductively for $1\leq\nu\leq n-1$
the common zero-set $Y_\nu$ of $s_1,\cdots,s_\nu$ with the structure
sheaf
$$
{\mathcal O}_{Y_\nu}:=\left({\mathcal O}_Y\left/\sum_{j=1}^\nu
s_j{\mathcal O}_Y\right)\right|_{Y_\nu}
$$
is regular of complex codimension $\nu$ in $Y$ and we end up with the inequality
$$
\displaylines{\dim_{\mathbb C}\Gamma\left(Y,\,{\mathcal
I}_{p\varphi}\left(pL+nA\right)\right)\cr\geq
\dim_{\mathbb C}\Gamma\left(Y_{n-1},\,{\mathcal
I}_{\left(p\varphi|Y_{n-1}\right)}\left(\left(pL+nA\right)|_{Y_{n-1}}\right)\right).\cr}
$$
For this step we have to exclude $P_0$ from a subvariety $Z_{n-2}$
of dimension $\leq n-2$ in $Y$, because we have to exclude a finite
set in each $Y_1$ which would come together as the hypersurface
$Y_1$ varies to form a subvariety $Z_1$ of dimension $\leq 1$ in $Y$
(as one can argue with the quotients of coherent ideal sheaves by
non zero-divisors and with the primary decompositions for coherent
ideal sheaves). Likewise we have a subvariety $Z_k$ of dimension
$\leq k$ in $Y$ so that $Z_k$ intersects $Y_k$ in a finite number of
points and finally we have end up with a subvariety $Z_{n-2}$ of
dimension $\leq n-2$ in $Y$ which intersects $Y_{n-2}$ in a finite
number of points and we impose the condition that $P_0$ does not
belong to $Z_{n-2}$.

\medbreak Since $Y_{n-1}$ is a curve, all coherent ideal sheaves on
it are principal and are locally free and they come from holomorphic
line bundles. We can choose $s_1,\cdots,s_{n-1}$ so generically that
$Y_{n-1}$ is disjoint from $Z_{n-2}$.  For this step we need to make sure
that $P_0$ does not belong to $Z_{n-2}$.

\medbreak\noindent(5.2.2) {\it Application of the Theorem of
Riemann-Roch to a Curve and Comparing Contributions from the Curvature
Current and the Multiplier Ideal Sheaf.} Let $c$ be the nonnegative number
$$\int_{Y_{n-1}}R=\int_Y
R\wedge\left(\omega_A\right)^{n-1},$$ where $\omega_A$ is the curvature form of some smooth positively curved metric of the ample line bundle $A$.  Then
$$
\displaylines{(5.2.2.1)\qquad\qquad\dim_{\mathbb
C}\Gamma\left(Y,\,{\mathcal
I}_{p\varphi}\left(pL+nA\right)\right)\hfill\cr\geq
\dim_{\mathbb C}\Gamma\left(Y_{n-1},\,{\mathcal
I}_{\left(\left.p\varphi\right|Y_{n-1}\right)}\left(\left(pL+nA\right)|_{Y_{n-1}}\right)\right)\cr
\geq 1-{\rm genus}\left(Y_{n-1}\right)+
nA\cdot Y_{n-1}
\cr+\sum_{j=1}^J\left(p\tau_j-\left\lfloor
p\tau_j\right\rfloor\right) V_j\cdot A^{n-1}+p\int_{Y_{n-1}}R,\cr}
$$
where the last identity is from the theorem of Riemann-Roch applied
to the regular curve $Y_{n-1}$ and the locally free sheaf
$${\mathcal
I}_{\left(\left.p\varphi\right|Y_{n-1}\right)}\left(\left(pL+nA\right)|_{Y_{n-1}}\right)\leqno{(5.2.2.2)}$$
on $Y_{n-1}$.  Note that, though for each $p$, there is a line bundle $E_p$ on $Y_{n-1}$ associated to the locally free sheaf (5.2.2.2) the line bundle $E_p$ in general is not the $p$-th tensor power of some fixed line bundle $E$ independent of $p$ even just for $p$ sufficiently large or an infinite sequence of distinct positive integers $p$.

\medbreak When $R$ is not identically zero, the nonnegative number $c$ is strictly positive and we conclude from (5.2.2.1) that
$$\dim_{\mathbb
C}\Gamma\left(Y,\,{\mathcal
I}_{p\varphi}\left(pL+nA\right)\right)\geq  1-{\rm genus}\left(Y_{n-1}\right)+
nA\cdot Y_{n-1}+pc$$
goes to $\infty$ as $p\to\infty$.

\medbreak Now assume the other case where $R=0$ and $J=\infty$.   By using arguments of diophantine approximation as explained below, we conclude that there exists some sequence $\left\{p_\nu\right\}_{\nu\in{\mathbb N}}$ of positive integers such that
$$\sum_{j=1}^J\left(p_\nu\tau_j-\left\lfloor
p_\nu\tau_j\right\rfloor\right) \to\infty\quad{\rm as}\ \ \nu\to\infty.\leqno{(5.2.2.3)}
$$
From (5.2.2.1) it follows that
$$
\lim_{\nu\to\infty}\dim_{\mathbb C}\Gamma\left(Y, {\mathcal I}_{p_\nu\varphi}\left(p_\nu L+n A\right)\right)=\infty.
$$
This ends the proof of Proposition (5.2) after we explain the derivation of (5.2.2.3), which we now do.  When all $\gamma_j$ are rational, the statement (5.2.2.3) follows from the convergence of $\sum_{j=1}^J\gamma_j$ and a simple comparison test for divergence or convergence of a series of positive terms.  When at least one of $\gamma_j$ is irrational, the statement (5.2.2.1) follows from the Corollary (5.2.2.4) to Kronecker's diophantine approximation which is a consequence of Kronecker's theorem listed below as (5.2.2.3) and given as Theorem 444 on p.382 of  [Hardy-Wright 1960].   A derivation of Corollary (5.2.2.4) from Kronecker's diophantine approximation (5.2.2.3) can be found, for example, in [Siu 2006, \S5].

\medbreak\noindent(5.2.2.3) {\it Theorem (Kronecker).}  Let
$a_1,\cdots,a_N$ be ${\mathbb Q}$-linearly independent real numbers.
Let $b_1,\cdots,b_N\in{\mathbb R}$.  Let $\varepsilon, T$ be
positive numbers. Then we can find $t>T$ and integers
$x_1,\cdots,x_N$ such that
$\left|ta_j-b_j-x_j\right|\leq\varepsilon$ for $1\leq j\leq N$.

\medbreak\noindent(5.2.2.4) {\it Corollary to Kronecker's Diophantine Approximation.}  Let $\gamma_j$ ($1\leq
j<\infty$) be a sequence of positive numbers and $\Lambda$ be a
positive integer such that $1,\gamma_1,\cdots,\gamma_\Lambda$ are
${\mathbb Q}$-linearly independent and
$$
\gamma_j=\sum_{\lambda=1}^{\Lambda}c_{j,\lambda}\gamma_\lambda
$$
for $\Lambda<j<\infty$, where $c_{j,k}\in{\mathbb Q}$.  For any
positive integer $N$ there exists some positive integer $m$ such
that $$m \gamma_j- \left\lfloor m
\gamma_j\right\rfloor\geq\frac{1}{4}\quad{\rm for\ \ }1\leq j\leq
N.$$

\bigbreak\noindent(5.3) {\it Remark.}  One important point about Proposition (5.2) is that the space
$\Gamma\left(Y, {\mathcal I}_{p_\nu\varphi}\left(p_\nu L+n A\right)\right)$ is used instead of
the space $\Gamma\left(Y, p_\nu L+n A\right)$.  For the application of Proposition (5.2) in the analytic proof of the finite generation of the canonical ring, the complex manifold $Y$ will be a nonsingular hypersurface of the compact complex algebraic manifold $X$ of general type.  The conclusion about the dimension of $\Gamma\left(Y, {\mathcal I}_{p_\nu\varphi}\left(p_\nu L+n A\right)\right)$ is used in producing nonzero elements of $\Gamma\left(Y, {\mathcal I}_{p_\nu\varphi}\left(p_\nu L+n A\right)\right)$ whose extensions to $X$ can be chosen to vanish to high order at a prescribed point of $Y$.   The high-order roots of the absolute-value-squares of such extensions are used to construct metrics for $K_X$ with nonnegative curvature current and extra singularities at a prescribed point of $Y$ but no more than the singularities of the metric of minimum singularities at a generic point of $Y$ plus a small prescribed error.  Such metrics hold the key to the proof of the precise achievement of the stable vanishing order across $Y$.  Only when we use the space
$\Gamma\left(Y, {\mathcal I}_{p_\nu\varphi}\left(p_\nu L+n A\right)\right)$ instead of
the space $\Gamma\left(Y, p_\nu L+n A\right)$, can we achieve the purpose of extending elements of $\Gamma\left(Y, {\mathcal I}_{p_\nu\varphi}\left(p_\nu L+n A\right)\right)$ to all of $X$ by using the vanishing theorem of Kawamata-Viehweg-Nadel for multiplier ideal sheaves.  The general type property of $X$ is needed to handle the small fraction of $nA$ in the process, just as in the proof of the deformational invariance of the plurigenera for the case of general type [Siu 1998].

\bigbreak\noindent(5.4) {\it Proposition (Extension of Modified Pluricanonical Section from Hypersurface of First Case of Dichotomy after Fixed Ample Twisting).}  Let $X$ be a compact complex algebraic manifold of general type of complex dimension $n$ and $e^{-\varphi}=\frac{1}{\Phi}$ be the metric of minimum singularity with curvature current $\Theta=\Theta_\varphi$.  Let $Y$ be a nonsingular hypersurface of $X$ and $\alpha$ be the generic stable vanishing order across $Y$.   Assume that the modified restriction $\Theta_\varphi-\alpha\left[Y\right]$ as defined in (3.3) belongs to the first case of the dichotomy in the sense that in the canonical decomposition
$$\Theta_\varphi-\alpha\left[Y\right]=\sum_{j=1}^J\gamma_j\left[V_j\right]+R,$$
with $\gamma_j>0$ and $V_j$ being a hypersurface in $Y$ and the Lelong number of $R$ vanishing outside a countable union of subvarieties of codimension $\geq 2$ in $Y$, either $J$ is infinite or $R$ is nonzero. 
Then there exists an ample line bundle $A_0$ on $X$ (which depends only on $X$ and $Y$) with the following property.  For any ample line bundle $A$ with $A-A_0$ ample, the complex dimension of
$$
\Gamma\left(X, {\mathcal I}_{p_\nu\left(\varphi-\alpha\log\left|s_Y\right|^2\right)}\left(p_\nu\left(K_X-\alpha Y\right)+A\right)\right)\Big|_Y 
$$
goes to $\infty$ as $\nu\to\infty$ for some increasing sequence $\left\{p_\nu\right\}_{\nu\in{\mathbb N}}$
of positive integers.

\medbreak\noindent{\it Proof.}  Let $A_1$ be a very ample line bundle on $X$ such that  $A_1-\left(K_X-\alpha Y\right)$ is ample on $X$ and $\left(A_1\Big|_Y\right)-K_Y$ is very ample on $Y$.  Let $\psi$ be the potential for the current $\Theta-\alpha\left[Y\right]$ which is the modified restriction of $\Theta$ to $Y$ in the sense described in (3.3) so that $\Theta-\alpha\left[Y\right]=\frac{\sqrt{-1}}{2\pi}\,\partial\bar\partial\psi$.   Then by (5.2) there exists some sequence $\left\{p_\nu\right\}_{\nu\in{\mathbb N}}$ of positive integers such that
$$
\lim_{\nu\to\infty}\dim_{\mathbb C}\Gamma\left(Y, {\mathcal I}_{p_\nu\varphi}\left(p_\nu L+(n-1) A_1\right)\right)=\infty.\leqno{(5.4.1)}
$$

\medbreak We choose an integer $\hat p>n$ such that $(3+\hat p-n)A_1$, $\tilde p A_1-Y$ and $\tilde p A_1$ are all globally free on $Y$.  Let  $\hat p=\left(1+\left\lceil\alpha\right\rceil-\alpha\right)\tilde p$.  We can construct multi-valued holomorphic sections $\sigma_1,\cdots,\sigma_\ell$ of $\hat p A_1$ such that their common zero-set is $Y$ and the vanishing order of 
$\sum_{j=1}^\ell\left|\sigma_j\right|^2$ across $Y$ is precisely $2\left(1+\alpha\right)$ at every point of $Y$ by setting $\sigma_j$ to be $\left(\hat\sigma s_Y\right)^{1+\alpha}\left(\tilde\sigma\right)^{\left\lceil\alpha\right\rceil-\alpha}$ as $\hat\sigma$ runs through a finite ${\mathbb C}$-basis of $\Gamma\left(X, \tilde p A_1-Y\right)$ and $\tilde\sigma$ runs through a finite ${\mathbb C}$-basis of $\Gamma\left(X, \tilde pA_1\right)$.  We now introduce the metric
$$
\frac{h_{A_1}\left|s_Y\right|^{2(p+1)\alpha}}{\Phi^p\sum_{j=1}^\ell\left|\sigma_j\right|^2}
$$
of $p\left(K_X-\alpha Y\right)-\alpha Y+\left(1+\hat p\right)A_1$ whose multiplier ideal sheaf is equal to ${\mathcal I}_{p\left(\varphi-\alpha\log\left|s_Y\right|^2\right)}\left({\mathcal I}d_Y\right)$, where $h_{A_1}$ is a smooth strictly positive curved metric of $A_1$. By the vanishing theorem of Kawamata-Viehweg-Nadel,
$$
H^1\left(X, {\mathcal I}_{p\left(\varphi-\alpha\log\left|s_Y\right|^2\right)}\left({\mathcal I}d_Y\right)\left(p\left(K_X-\alpha Y\right)+\left(2+\hat p\right)A_1\right)\right)=0,\leqno{(5.4.2)}
$$
because $A_1-K_Y$ is ample on $Y$.  An element of ${\mathcal I}_{p\psi}$ on $Y$ can be naturally regarded as an element of $${\mathcal I}_{p\left(\varphi-\alpha\log\left|s_Y\right|^2\right)}\left/\left(
{\mathcal I}_{p\left(\varphi-\alpha\log\left|s_Y\right|^2\right)}\left({\mathcal I}d_Y\right)\right)\right.,$$
for example, by using the extension theorem of Ohsawa-Takegoshi [Ohsawa-Takegoshi 1987].  From (5.4.2) it follows that the restriction map
$$
\displaylines{\Gamma\left(X, {\mathcal I}_{p\left(\varphi-\alpha\log\left|s_Y\right|^2\right)}
\left(p\left(K_X-\alpha Y\right)+\left(2+\hat p\right)A_1\right)\right)\cr\to
\Gamma\left(X, \left({\mathcal I}_{p\left(\varphi-\alpha\log\left|s_Y\right|^2\right)}\left/\left(
{\mathcal I}_{p\left(\varphi-\alpha\log\left|s_Y\right|^2\right)}\left({\mathcal I}d_Y\right)\right)\right.\right)
\left(p\left(K_X-\alpha Y\right)+\left(2+\hat p\right)A_1\right)\right)\cr}
$$
is surjective.  Thus the the restriction map
$$
\displaylines{(5.4.3)\qquad\Gamma\left(X, {\mathcal I}_{p\left(\varphi-\alpha\log\left|s_Y\right|^2\right)}
\left(p\left(K_X-\alpha Y\right)+\left(2+\hat p\right)A_1\right)\right)\hfill\cr\hfill\to
\Gamma\left(Y, {\mathcal I}_{p\psi}\left(p\left(K_X-\alpha Y\right)+\left(2+\hat p\right)A_1\right)\right)\qquad\qquad\cr}
$$
is surjective.  Since $\left(3+\hat p-n\right)A_1$ is globally free on $Y$ and the product of an element of $\Gamma\left(Y, {\mathcal I}_{p\psi}\left(p\left(K_X-\alpha Y\right)+(n-1)A_1\right)\right)$ and an element of $\Gamma\left(Y, \left(3+\hat p-n\right)A_1\right)$ is an element of $\Gamma\left(Y, {\mathcal I}_{p\psi}\left(p\left(K_X-\alpha Y\right)+\left(2+\hat p\right)A_1\right)\right)$, it follows from (5.4.1) and the surjectivity of (5.4.3) that
$$
\lim_{\nu\to\infty}\dim_{\mathbb C}\left(\Gamma\left(X, {\mathcal I}_{p_\nu\left(\varphi-\alpha\log\left|s_Y\right|^2\right)}
\left(p_\nu\left(K_X-\alpha Y\right)+\left(2+\hat p\right)A_1\right)\right)\bigg|_Y\right)=\infty.
$$

\medbreak Let $A_2$ be a holomorphic line bundle on $X$ which is sufficient ample so that for any ample line bundle $E$ the line bundle $A_2+E$ is globally free.  For example, we can choose $A_2$ such that for every point $P$ of $X$ there exist elements of $\Gamma\left(X, A_2-K_X\right)$ which vanish to order $\geq n+1$ at $P$ and whose common zero-set in $X$ consists of just the single point $P$.  We now set $A_0$ to be any ample line bundle over $X$ such that $A_0-\left(2+\tilde p\right)A_1-A_2$ is ample.   Since for any holomorphic line bundle $A$ on $X$ with $A-A_0$ ample on $X$ the product of any element of $$\Gamma\left(X, {\mathcal I}_{p_\nu\left(\varphi-\alpha\log\left|s_Y\right|^2\right)}
\left(p_\nu\left(K_X-\alpha Y\right)+\left(2+\hat p\right)A_1\right)\right)$$ and any element of $\Gamma\left(X, A-\left(2+\tilde p\right)A_1\right)$ is an element of $$\Gamma\left(X, {\mathcal I}_{p_\nu\left(\varphi-\alpha\log\left|s_Y\right|^2\right)}
\left(p_\nu\left(K_X-\alpha Y\right)+A\right)\right),$$ it follows from the very ample property of $A-\left(2+\tilde p\right)A_1$ that the complex dimension of
$$
\Gamma\left(X, {\mathcal I}_{p_\nu\left(\varphi-\alpha\log\left|s_Y\right|^2\right)}\left(p_\nu\left(K_X-\alpha Y\right)+A\right)\right)\Big|_Y 
$$
goes to $\infty$ as $\nu\to\infty$.  Q.E.D.

\medbreak\noindent(5.4.1) {\it Remark.}  In the proof of (5.4) we used the metric
$$
\frac{h_{A_1}\left|s_Y\right|^{2(p+1)\alpha}}{\Phi^p\sum_{j=1}^\ell\left|\sigma_j\right|^2}
$$
of $p\left(K_X-\alpha Y\right)-\alpha Y+\left(1+\hat p\right)A_1$ to produce the multiplier ideal sheaf ${\mathcal I}_{p\left(\varphi-\alpha\log\left|s_Y\right|^2\right)}\left({\mathcal I}d_Y\right)$, which contains precisely the additional factor $\left({\mathcal I}d_Y\right)$, compared to the multiplier ideal sheaf ${\mathcal I}_{p\left(\varphi-\alpha\log\left|s_Y\right|^2\right)}$ of the metric $e^{-p\varphi}\left|s_Y\right|^{2p\alpha}$.  The reason why this can be done is the use of the sufficiently ample line bundle $A_1$ over $X$.  

\medbreak Of course, later we are able to take care of the effect of the use of $A_1$, because the dimension of the space of sections under consideration goes to infinity as $p_\nu\to\infty$ so that we can divide the coefficient of $A_1$ by an arbitrarily large number in the proof of (5.5) below.   In the second case of the dichotomy to be discussed in \S6 there will be no such way to get rid of the effect of $A_1$ due to the uniform bound of the dimension of the space of sections under consideration.  

\medbreak For the application of the vanishing theorem of Kamawata-Viehweg-Nadel to extend sections, in \S6 below we will be forced to produce an additional factor ${\mathcal J}$ of the multiplier ideal sheaf which is only equal to ${\mathcal I}d_Y$ at a generic point of $Y$ instead of everywhere.  This is obtained by using a metric $\frac{\left|s_Y\right|^{2p\alpha}}{\tilde\Phi^p}$ not involving any ample line bundle over $X$, where $\tilde\Phi$ is a finite sum of absolute-value-squares of multi-valued holomorphic sections of $K_X$ with appropriately higher vanishing orders across $Y$ than the $m$-th roots of elements of a basis of the space of $m$-canonical sections over $X$.  The fact that ${\mathcal J}$ is only equal to ${\mathcal I}d_Y$ at a generic point of $Y$ instead of everywhere poses the difficulty of ``additional vanishing''.  This makes the second case of the dichotomy more complicated than the first case of the dichotomy and will be dealt with in \S6.

\bigbreak\noindent(5.5) {\it Proposition (Metric of Additional Singularity on Hypersurface of First Case of Dichotomy after Fixed Ample Twisting).}  Let $X$ be a compact complex algebraic manifold of general type of complex dimension $n$ and $e^{-\varphi}=\frac{1}{\Phi}$ be the metric of minimum singularity with curvature current $\Theta=\Theta_\varphi$.  Let $Y$ be a nonsingular hypersurface of $X$ and $\alpha$ be the generic stable vanishing order across $Y$.   Assume that the modified restriction $\Theta_\varphi-\alpha\left[Y\right]$ as defined in (3.3) belongs to the first case of the dichotomy in the sense that in the canonical decomposition
$$\Theta_\varphi-\alpha\left[Y\right]=\sum_{j=1}^J\gamma_j\left[V_j\right]+R,$$
with $\gamma_j>0$ and $V_j$ being a hypersurface in $Y$ and the Lelong number of $R$ vanishing outside a countable union of subvarieties of codimension $\geq 2$ in $Y$, either $J$ is infinite or $R$ is nonzero. 
Then for every point $P_0$ and any $\varepsilon>0$ and any $N\in{\mathbb N}$ there exists a metric $$e^{-\varphi_{Y,P_0,\varepsilon,N}}=\frac{1}{\Phi_{Y,P_0,\varepsilon,N}}$$ of $\hat p K_X$ for some $\hat p=\hat p_{Y,P_0,\varepsilon,N}\in{\mathbb N}$ such that
\begin{itemize}
\item[(i)] $\Phi_{Y,P_0,\varepsilon,N}$ is a finite sum of absolute-value-squares of multi-valued holomorphic sections of $\hat p K_X$ over $X$,
\item[(ii)] the vanishing order of $\Phi_{Y,P_0,\varepsilon,N}$ across $Y$ at some generic point of $Y$ is a number in the interval $\left[2\hat p\alpha,2\hat p\alpha+\varepsilon\right)$ ,
\item[(iii)] the vanishing order of $\Phi_{Y,P_0,\varepsilon,N}$ on $X$ at $P_0$ is at least $2N$. 
\item[(iv)] the curvature current of the metric $e^{-\varphi_{Y,P_0,\varepsilon,N}}$ is strictly positive on $X$.
\end{itemize}

\medbreak\noindent{\it Proof.}  Since $X$ is of general type, we have $K_X=D+B$ for some effective ${\mathbb Q}$-divisor $D$ and some ample ${\mathbb Q}$-line bundle $B$.  Let $D=aY+E$ with $a\geq 0$ and $Y$ not contained in the support of $E$.  Then $K_X=aY+E+B$.   Since $\alpha$ is the generic stable vanishing order across $Y$, it follows that $a\geq\alpha$. 

\medbreak By Proposition (5.4) on the extension of modified pluricanonical sections from a hypersurface of the first case of dichotomy after fixed ample twisting, we have an ample line bundle $A_0$ on $X$ such that, for any ample line bundle $A$ with $A-A_0$ ample, the complex dimension of
$$
\Gamma\left(X, {\mathcal I}_{p_\nu\left(\varphi-\alpha\log\left|s_Y\right|^2\right)}\left(p_\nu\left(K_X-\alpha Y\right)+A\right)\right)\Big|_Y 
$$
goes to $\infty$ as $\nu\to\infty$ for some increasing sequence $\left\{p_\nu\right\}_{\nu\in{\mathbb N}}$
of positive integers.  Thus, for any $q\in{\mathbb N}$ there exists some $\tilde p_q$ such that we can find an element
$$
\sigma_q\in \Gamma\left(X, {\mathcal I}_{\tilde p_q\left(\varphi-\alpha\log\left|s_Y\right|^2\right)}\left(\tilde p_q\left(K_X-\alpha Y\right)+2A_0\right)\right)
$$
whose restriction to $Y$ is not identically zero and which vanishes to order $\geq q$ on $X$ at $P_0$.

\medbreak We choose $\ell\in{\mathbb N}$ with $\frac{\alpha-a}{\ell}<\frac{\varepsilon}{4}$.  Choose $\hat\ell\in{\mathbb N}$ such that $\frac{1}{\ell}B-\frac{2}{\hat\ell}A_0$ is an ample ${\mathbb Q}$-line bundle over $X$.  Choose $q\in{\mathbb N}$ with $q>\hat\ell N$.  Choose a finite number of multi-valued holomorphic sections $\tau_1,\cdots,\tau_k$ of the ${\mathbb Q}$-line bundle $\frac{1}{\ell}B-\frac{2}{\hat\ell}A_0$ over $X$ such that the curvature current of the metric 
$$
\frac{1}{\sum_{j=1}^k\left|\tau_j\right|^2}
$$
of the ${\mathbb Q}$-line bundle $\frac{1}{\ell}B-\frac{2}{\hat\ell}A_0$ is strictly positive on $X$.   Choose $N_0\in{\mathbb N}$ such that the generic vanishing order of the $N_0$-truncation $\Phi_{N_0}$ of $\Phi$ is a number $2\gamma$ in the interval $\left[2\alpha,2\alpha+\frac{\varepsilon}{2}\right)$.  Let $\hat p_{Y,P_0,\varepsilon,N}\in{\mathbb N}=1+\tilde p_q$.  Then
$$
\Phi_{Y,P_0,\varepsilon,N}=\left|\left(\sigma_q\right)^{\frac{1}{\hat\ell}}
\left(s_E\right)^{\frac{1}{\ell}}\left(s_Y\right)^{\frac{\alpha-a}{\ell}}\right|^2\left(\sum_{j=1}^k\left|\tau_j\right|^2\right)\left(\Phi_{N_0}\right)^{1-\frac{1}{\ell}}\left|s_Y\right|^{2\hat p\alpha}
$$
satisfies the requirement, because the generic vanishing order of
$\Phi_{Y,P_0,\varepsilon,N}$ across $Y$ is $2\left(\hat p\alpha+\gamma+\frac{\alpha-a}{\ell}\right)$ and $\left(\gamma-\alpha\right)+\frac{\alpha-a}{\ell}<\frac{\varepsilon}{2}$.
Q.E.D.

\bigbreak\noindent(5.5.1) {\it Remark on Holomorphic Dependence.}  In the proof of (5.5), instead of using a single
$$
\sigma_q\in \Gamma\left(X, {\mathcal I}_{\tilde p_q\left(\varphi-\alpha\log\left|s_Y\right|^2\right)}\left(\tilde p_q\left(K_X-\alpha Y\right)+2A_0\right)\right)
$$
we could also have used a maximal ${\mathbb C}$-linear independent subset $\hat\sigma_{q,1},\cdots,\hat\sigma_{q,t_q}$ of 
$$
\Gamma\left(X, \left(\left({\mathfrak m}_{X,P_0}\right)^q\cap{\mathcal I}_{\tilde p_q\left(\varphi-\alpha\log\left|s_Y\right|^2\right)}\right)\left(\tilde p_q\left(K_X-\alpha Y\right)+2A_0\right)\right)\Big|_Y,
$$
whose restrictions to $Y$ are still ${\mathbb C}$-linearly independent, to form the metric
$$
\hat\Phi_{Y,P_0,\varepsilon,N}=\left(\sum_{j=1}^{t_q}
\left|\hat\sigma_{q,j}\right|^2\right)^{\frac{1}{\hat\ell}}
\left|\left(s_E\right)^{\frac{1}{\ell}}\left(s_Y\right)^{\frac{\alpha-a}{\ell}}\right|^2\left(\sum_{j=1}^k\left|\tau_j\right|^2\right)\left(\Phi_{N_0}\right)^{1-\frac{1}{\ell}}\left|s_Y\right|^{2\hat p\alpha}
$$
of $\hat p K_X$.
The advantage of using a maximal ${\mathbb C}$-linear independent set to construct $
\hat\Phi_{Y,P_0,\varepsilon,N}$ instead of using $\Phi_{Y,P_0,\varepsilon,N}$ is that when $P_0$ varies in some appropriate Zariski open subset of $Y$, the Lelong set $E_c$ (for $\eta_1 N\leq c\leq\eta_2 N$ for suitable $0<\eta_1<\eta_2<1$ and for $N$ sufficiently large) of the curvature current of the metric
$$
\left.\frac{\left|s_Y\right|^{2\left(\hat p\alpha+\gamma+\frac{\alpha-a}{\ell}\right)}}{\hat\Phi_{Y,P_0,\varepsilon,N}}\right|_Y
$$
of the line bundle $\hat K_X-\left(\hat p\alpha+\gamma+\frac{\alpha-a}{\ell}\right)Y$ on $Y$ varies holomorphically as a function of $P_0$.  (See (3.1) for the definition of the Lelong set $E_c$.)

\bigbreak\noindent(5.5.2) {\it Remark on the Inevitability of Eventual Occurrence of the Second Case of the Dichotomy.}  In the arguments from (5.2) to (5.5) our goal is to take a step toward producing holomorphic pluricanonical sections $s\in\Gamma\left(X, mK_X\right)$ on $X$ to achieve the generic stable vanishing order $\alpha$ across the hypersurface $Y$ in $X$ by assuming that $Y$ is in the first case of the dichotomy.  That is, the vanishing order of $s$ across $Y$ is $m\alpha$ at a generic point of $Y$.  
The step we take is to produce a positively curved metric of $K_X$ whose reciprocal has high-order vanishing at a prescribed point of $Y$ and generic vanishing order across $Y$ as close to $2\alpha$ as prescribed.

\medbreak However, the construction of such a metric is as good as we can go with this method.  We can replace $Y$ by a subvariety $Y^\prime$ in $Y$ defined by high-order vanishing of the reciprocal of such a metric and apply the argument to $Y^\prime$ instead of to $Y$ and get another metric if $Y_1$ belongs to the first case of the dichotomy in a naturally corresponding sense.    We can keep repeating this argument of constructing metrics for $Y,\,Y^\prime,\,Y^{\prime\prime},\cdots$ (as long as each of $Y^\prime,\,Y^{\prime\prime}, \cdots$ belongs to the first case of the dichotomy), but, without some additional different techniques, we cannot use such arguments to finally get  actual holomorphic pluricanonical sections $s\in\Gamma\left(X, mK_X\right)$ achieving the generic stable vanishing order across $Y$.   The reason is as follows.

\medbreak The current arguments use extension to $X$ of sections on $Y$ by the vanishing theorem of Kawamata-Viehweg-Nadel.   The vanishing theorem of Kawamata-Viehweg-Nadel takes a line bundle $E$ on $X$ with a metric whose curvature current is strictly positive and gives the vanishing of cohomology of positive degree with the multiplier ideal sheaf ${\mathcal I}$ as the coefficient for the line bundle $E+K_X$.  The key point is that a copy of $K_X$ is added to $E$ to give the vanishing of the higher-degree cohomology $H^p\left(X, {\mathcal I}\left(E+K_X\right)\right)$ for $p\geq 1$.   This is applied to $E=m\left(K_X-\alpha Y\right)$.  

\medbreak To use $H^1\left(X, {\mathcal I}\left(m\left(K_X-\alpha Y\right)+K_X\right)\right)=0$ to extend a nontrivial section $\sigma$ on $Y$, the section $\sigma$ should be an element of $\Gamma\left(Y, m\left(K_X-\alpha Y\right)+K_X\right)$ which locally belongs to an appropriate ideal sheaf.   The result $\tilde\sigma$ of the extension would be an element of  $\Gamma\left(X, m\left(K_X-\alpha Y\right)+K_X\right)$.  The element in $\Gamma\left(X, (m+1)K_X\right)$ defined by $\tilde\sigma$ would have a generic vanishing order across $Y$ equal to $m\alpha$, which is less than the smallest possible number $(m+1)\alpha$, giving a contradiction, unless $\alpha=0$.

\medbreak The best we can get with such arguments is only the construction of a positively curved metric of $K_X$ whose reciprocal has high-order vanishing at a prescribed point of $Y^{(k)}$ and generic vanishing orders lexicographically across $Y, Y^\prime, Y^{\prime\prime},\cdots,Y^{(k-1)}$ as close to the expected smallest value as prescribed (as long as each of $Y^\prime,\,Y^{\prime\prime}, \cdots, Y^{(k-1)}$ belongs to the first case of the dichotomy).   Finally, we eventually get to $Y^{(k)}$ which belongs to the second case of the dichotomy, in the lexicographical sense from the nested sequence
$$Y^{(k)}\subset Y^{(k-1)}\subset\cdots\subset Y^{\prime\prime}\subset Y^\prime\subset Y\subset X.$$
It means that the curvature current on $Y$ is modified from the curvature $\Theta$ of $\frac{1}{\Phi}$ by taking away a rational multiple of $Y$ and restrict to $Y$ and then is modified by taking away a rational multiple of $Y^\prime$ and restrict to $Y^\prime$, et cetera, until the final modification obtained by taking away a rational multiple of $Y^{(k-1)}$ and restrict to $Y^{(k)}$.  Because of the singularities, a rigorous treatment would require blowing up to regular hypersurfaces in normal crossing at each stage.  This meaning of the second case of the dichotomy by inductive lexicographical description was already sketchily mentioned in [Siu 2008, (2.2)]. 

\bigbreak\noindent(5.6) {\it Holomorphic Family of Embedded Subvarieties of High Stable Vanishing Order Inside a Hypersurface Whose Stable Vanishing Order Is Not Yet Known to Be Achieved.}   We continue the discussion in (5.5) and use the same assumptions and notations as in (5.5.2).  By the discussion in (5.5.1) the nested subvarieties 
$$Y^{(k)}\subset Y^{(k-1)}\subset\cdots\subset Y^{\prime\prime}\subset Y^\prime\subset Y\subset X$$
defined by the metrics actually form holomorphic families outside appropriate Zariski open subsets.  From the discussion in (4.5.1) applied to the fibrations inside each $Y^\prime, Y^{\prime\prime},\cdots,Y^{(k-1)}$ and (4.6), we can decrease the discrepancy subspace until we get to a finite union $Z$ of subvarieties which belong to the second case of the dichotomy.  The argument to decrease the discrepancy subspace comes from the use of a sufficiently ample line bundle from the base of the fibration.  As a result, the difficulty of  ``additional vanishing'' described in (6.1) below from the use of pluricanonical sections with appropriate generic vanishing orders across specified subvarieties does not occur outside $Z$.

\bigbreak\noindent{\sc\S6.  Subspaces of Minimum Additional Vanishing for the Second Case of the Dichotomy}

\bigbreak\noindent(6.1) {\it Difficulty of Additional Vanishing}.  Let $X$ be a compact complex algebraic manifold of complex dimension $n$ of general type.   Suppose the common zero-set of all $m$-canonical sections of $X$ for all $m\in{\mathbb N}$ is $Z=\bigcup_{j=1}^\ell\alpha_j Y_j$ for some collection $\left\{Y_j\right\}_{j=1}^\ell$ are nonsingular hypersurfaces of $X$ in normal crossing.  Let $\alpha_j$ be the generic stable vanishing order for $Y_j$ for $1\leq j\leq\ell$.  

\medbreak By relabelling $Y_\nu$ for $1\leq\nu\leq\ell$, we assume that there exists $1<\ell^\prime\leq\ell$ such that
\begin{itemize}\item[(i)] for $1\leq\nu\leq\ell^\prime$ there does not exist any multi-valued holomorphic section of $K_X$ over $X$ whose generic vanishing order across $Y_\nu$ is $\alpha_\nu$, and
\item[(ii)] for $\ell^\prime<\nu\leq\ell$ there exists some multi-valued holomorphic section $s_\nu$ of $K_X$ over $X$ whose generic vanishing order across $Y_\nu$ is precisely $\alpha_\nu$.
\end{itemize}
We let $s_{\ell+1},\cdots,s_{\tilde\ell}$ be multi-valued holomorphic sections of $K_X$ whose common zero-set is $Z$.  Choose $0<\tilde\delta<1$ such that the generic common vanishing order of $s_{\ell+1},\cdots,s_{\tilde\ell}$ across  $Y_\nu$ is $>\alpha_\nu+\tilde\delta$ for $1\leq\nu\leq\ell^\prime$.
We can take a multi-valued holomorphic section $s_\nu$ of $K_X$ over $X$ such that the generic vanishing order of $s_\nu$ across $Y_\nu$ is $\alpha_\nu+\delta_\nu$ with $0<\delta_\nu<\hat\delta$ for $1\leq\nu\leq\ell^\prime$.    The set of  the multi-valued holomorphic sections $s_1,\cdots,s_{\tilde\ell} $ of $K_X$ over $X$ satisfies the following conditions.
\begin{itemize}\item[(i)]
The common zero-set of the multi-valued holomorphic sections $s_1,\cdots,s_{\tilde\ell} $ of $K_X$ over $X$ is $Z$. \item[(ii)] The common generic vanishing order of $s_1,\cdots,s_{\tilde\ell}$ across $Y_\nu$ is $\alpha_\nu+\delta_\nu$ for $1\leq\nu\leq\ell^\prime$.  \item[(iii)]  The common generic vanishing order of $s_1,\cdots,s_{\tilde\ell}$ across $Y_\nu$ is $\alpha_\nu$ for $\ell^\prime<\nu\leq\ell$.\end{itemize}
For any $Y_\nu$ (with $1\leq\nu\leq\ell$) in the first case of the dichotomy, we use the techniques in \S5 to handle so that eventually we end up with the second case of the dichotomy.  We now discuss how to handle the second case of the dichotomy.  

\medbreak We assume that each $Y_\nu$ for $1\leq\nu\leq\ell^\prime$ belongs to the second of the dichotomy.  Moreover, we assume that  the common vanishing order of $s_1,\cdots,s_{\tilde\ell}$ across $Y_\nu$ at every point of $Y_\nu-\bigcup_{j=1}^{\ell^\prime}Y_j$ is $\alpha_\nu$ for $\ell^\prime<\nu\leq\ell$.  So for $1\leq\nu\leq\ell^\prime$ and the modified restriction of the curvature current $\Theta-\alpha_\nu Y_\nu$ is of the form $\sum_{j=1}^{J_\nu}\gamma_{\nu,j}\left[V_{\nu,j}\right]$ on $Y_\nu$.    It means that on $Y_\alpha$ the ${\mathbb Q}$-line bundle $K_X-\alpha_\nu Y_\nu$ is equal to a flat bundle $F$ on $Y_\alpha$ plus the line bundle $\sum_{j=1}^{J_\nu}\gamma_{\nu,j}\left[V_{\nu,j}\right]$ on $Y_\nu$.  Note that the second case of the dichotomy does not mean that the restriction of $K_X-\alpha_\nu Y_\nu$ to $Y_\nu$ is flat.  It only means that the curvature current on $Y_\nu$ of the restriction of $K_X-\alpha_\nu Y_\nu$ to $Y_\nu$ is the positive linear combination of only a finite number of hypersurfaces of $Y_\nu$.  There are two things that need to be done.
\begin{itemize}\item[(i)] We have to use Shokurov's technique to produce, from the special form of the curvature $\sum_{j=1}^{J_\nu}\gamma_{\nu,j}\left[V_{\nu,j}\right]$, a non identically zero holomorphic section $\rho_\nu$ of $q_\nu\left(K_X-\alpha_\nu Y_\nu\right)$ on $Y_\nu$ for some sufficiently large (and appropriately divisible) positive integer $q_\nu$.
\item[(ii)]  We have to extend  $\rho_\nu$ to an element of $\Gamma\left(X, q_\nu\left(K_X-\alpha_\nu Y_\nu\right)\right)$ by the vanishing theorem of Kawamata-Viehweg-Nadel.
\end{itemize}
Let us first handle Item (ii).    Item (i) will be discussed in (6.4) below.

\medbreak Since $X$ is of general type, we can write $K_X=D+B$, where $D$ is an effective ${\mathbb Q}$-divisor of $X$ and $B$ is an ample ${\mathbb Q}$-line bundle over $X$.  We can choose a multi-valued section $\sigma$ of $B$ over $X$ such that the coefficients of $Y_\nu$ in the divisor $D+{\rm div\,}\sigma$ are all distinct for $1\leq\nu\leq\ell^\prime$.   By replacing $s_\nu$ by $\left(s_\nu\right)^{\alpha_\nu}\left(\sigma\right)^{1-\alpha_\nu}$ by appropriately chosing rational numbers $0\leq\alpha_\nu<1$ for $1\leq\nu\leq\ell^\prime$, we can assume that the positive rational numbers $\frac{1+\alpha_\nu}{\delta_\nu}$ for $1\leq\nu\leq\ell^\prime$ are all distinct.
By relabelling $Y_\nu$ for $1\leq\nu\leq\ell^\prime$, we can assume without loss of generality that 
$$
\frac{1+\alpha_\nu}{\delta_\nu}>\frac{1+\alpha_\mu}{\delta_\mu}\ \ {\rm for}\ \ 1\leq\nu<\mu\leq\ell^\prime.
$$
Let $q_\nu$ be any positive integer such that $q_\nu-1>\frac{1+\alpha_\nu}{\delta_\nu}$ for $1\leq\nu\leq\ell^\prime$.
For $1\leq\nu\leq\ell^\prime$ we introduce the metric
$$
e^{-\varphi_\nu}=\frac{1}{\Phi^{q_\nu-1-\frac{1+\alpha_\nu}{\delta_\nu}}\left(\sum_{j=1}^{\tilde\ell}\left|s_j\right|^2\right)^{\frac{1+\alpha_\nu}{\delta_\nu}}}
$$
of $\left(q_\nu-1\right)K_X$ and the accompanying metric
$$
e^{-\tilde\varphi_\nu}=\frac{\left(\prod_{j=1}^\ell \left|s_{Y_j}\right|^{2a_j}\right)^{q_\nu}}{\Phi^{q_\nu-1-\frac{1+\alpha_\nu}{\delta_\nu}}\left(\sum_{j=1}^{\tilde\ell}\left|s_j\right|^2\right)^{\frac{1+\alpha_\nu}{\delta_\nu}}}
$$
of $$\left(q_\nu-1\right)\left(K_X-\sum_{j=1}^\ell a_j Y_j\right)-\sum_{j=1}^\ell a_j Y_j.$$
The support of the multiplier ideal sheaf $\widetilde{{\mathcal I}_\nu}$ of the metric $e^{-\tilde\varphi_\nu}$ of $$\left(q_\nu-1\right)\left(K_X-\sum_{j=1}^\ell a_j Y_j\right)-\sum_{j=1}^\ell a_j Y_j$$ is $\bigcup_{j=\nu}^{\ell^\prime}Y_j$ and outside the union of $\bigcup_{j=\nu+1}^{\ell^\prime}Y_j$ and a subvariety $F_\nu$ of codimension $\geq 1$ in $\bigcup_{j=\nu}^{\ell^\prime}Y_j$ the ideal sheaf $\widetilde{{\mathcal I}_\nu}$ is equal to ${\mathcal I}_{
\left(q_\nu-1\right)\varphi-q_\nu\sum_{j=1}^\ell\log\left|s_{Y_\nu}\right|^2}$ times the full ideal sheaf ${\mathcal I}d_{Y_\nu}$ of $Y_\nu$.    The main point is that in general $\widetilde{{\mathcal I}_\nu}$ is not equal to ${\mathcal I}_{
\left(q_\nu-1\right)\varphi-q_\nu\sum_{j=1}^\ell\log\left|s_{Y_\nu}\right|^2}\,{\mathcal I}d_{Y_\nu}$  and we can only conclude that there exists some nonnegative $e_\nu$ such that
$$
\left({\mathcal I}d_{F_\nu}\right)^{e_\nu}
{\mathcal I}_{
\left(q_\nu-1\right)\varphi-q_\nu\sum_{j=1}^\ell\log\left|s_{Y_\nu}\right|^2}\,{\mathcal I}d_{Y_\nu}
\subset \widetilde{{\mathcal I}_\nu}\subset{\mathcal I}_{
\left(q_\nu-1\right)\varphi-q_\nu\sum_{j=1}^\ell\log\left|s_{Y_\nu}\right|^2}\,{\mathcal I}d_{Y_\nu}.
$$
This is the phenomenon of ``additional vanishing'', which presents difficulties when we try to extend sections by using the vanishing theorem of Kawamata-Viehweg-Nadel in the standard way given in the next step. 

\medbreak For this next step we use the general type property of $X$ to slightly change the metric $e^{-\tilde\varphi_\nu}$ of $$\left(q_\nu-1\right)\left(K_X-\sum_{j=1}^\ell a_j Y_j\right)-\sum_{j=1}^\ell a_j Y_j$$ so that the curvature of the new metric is strictly positive while its multiplier ideal sheaf is not changed.  We then use the vanishing theorem of Kawamata-Viehweg-Nadel to get  the following vanishing of the cohomology
$$
\displaylines{H^1\left(X,\,\widetilde{{\mathcal I}_\nu}\left(q_\nu\left(K_X-\sum_{j=1}^\ell a_j Y_j\right)\right)\right)\cr=H^1\left(X,\,\widetilde{{\mathcal I}_\nu}\left(\left(\left(q_\nu-1\right)\left(K_X-\sum_{j=1}^\ell a_j Y_j\right)-\sum_{j=1}^\ell a_j Y_j\right)+K_X\right)\right)=0.}$$
The restriction map
$$\Gamma\left(X,\,q_\nu\left(K_X-\sum_{j=1}^\ell a_j Y_j\right)\right)\to
\Gamma\left(Y_\nu,\,\left({\mathcal O}_X\left/\widetilde{{\mathcal I}_\nu}\right.\right)\left(q_\nu\left(K_X-\sum_{j=1}^\ell a_j Y_j\right)\right)\right)$$
is surjective.  In general we cannot replace the element $$\rho_\nu\in\Gamma\left(Y_\nu,q_\nu\left(K_X-\sum_{j=1}^\ell a_j Y_j\right)\right)$$ by another element of $\Gamma\left(Y_\nu,q_\nu\-\left(K_X-\sum_{j=1}^\ell a_j Y_j\right)\right)$ whose local extension to $X$ belongs to $\left({\mathcal I}d_{F_\nu}\right)^{e_\nu}$ so that it can be naturally regarded as an element of
$$\Gamma\left(Y_\nu,\,\left({\mathcal O}_X\left/\widetilde{{\mathcal I}_\nu}\right.\right)\left(q_\nu\left(K_X-\sum_{j=1}^\ell a_j Y_j\right)\right)\right)$$
which can then be lifted to an element of
$$\Gamma\left(X,\,q_\nu\left(K_X-\sum_{j=1}^\ell a_j Y_j\right)\right).$$
We call this difficulty of the inability to extend $\rho_\nu$ from $Y_\nu$ to all of $X$ the difficulty of {\it additional vanishing}, because the additional vanishing orders of the multiplier ideal sheaf $\widetilde{{\mathcal I}_\nu}$ at points of $Y_\nu$ which are over and above that of the full ideal sheaf ${\mathcal I}d_{Y_\nu}$ on $X$.

\bigbreak\noindent(6.2) {\it  Technique of Restricting First to a Subspace and then Extending Only the Restriction.}   We now discuss one technique of handling the difficulty of ``additional vanishing'' for the second case of the dichotomy.  This technique is to consider the section $\hat\rho$ over $\bigcup_{j=1}^{\ell^\prime}Y_\nu$ formed by all $\rho_\nu$ ($1\leq\nu\leq\ell^\prime$) put together when $q_\nu$ is chosen to be equal to the same number $\hat q$ for $1\leq\nu\leq\ell^\prime$, to restrict $\hat\rho$ first to a subspace $W$ of $X$ whose support is contained in $\bigcup_{j=1}^{\ell^\prime}Y_\nu$ and then to extend to $X$ only $\hat\rho|_W$ {\it if $W$ turns out to be a subspace of $\bigcup_{j=1}^{\ell^\prime}Y_\nu$ instead of just a subspace of $X$} and the restriction of $\hat\rho$ to $W$ is nontrivial.

\medbreak In general the condition that $W$ is a subspace of $\bigcup_{j=1}^{\ell^\prime}Y_\nu$, instead of just a subspace of $X$, cannot be satisfied.   So we cannot directly use it.  We present it here simply as a motivation and as some background material to understand better the idea and the use of subspaces of ``minimum additional singularity'' which we will do in (6.3).   The idea of a subspace of ``minimum additional singularity''  is the same as doing the extension only from such a subspace.  Any section which we would like to extend has first to be restricted to such a subspace of ``minimum additional singularity'' and only the restriction is to be extended and not the original section.   Here in general we cannot directly use $W$ because we are constraining it to be inside the reduced subspace $\bigcup_{j=1}^{\ell^\prime}Y_\nu$.  In (6.3) we will use a subspace of ``minimum additional vanishing'' to be not necessarily inside the reduced subspace $\bigcup_{j=1}^{\ell^\prime}Y_\nu$. 

\medbreak Let $\hat q$ be any positive integer and $\gamma$ be a positive rational number such that $\alpha_\nu\left(\hat q-1\right)>\gamma\left(\alpha_\nu+\delta_\nu\right)$ for $1\leq\nu\leq\ell^\prime$.
We introduce the metric
$$
e^{-\check\varphi_\gamma}=\frac{1}{\Phi^{\hat q-1-\gamma}\left(\sum_{j=1}^{\tilde\ell}\left|s_j\right|^2\right)^\gamma}
$$
of $\left(\hat q-1\right)K_X$ and the accompanying metric
$$
e^{-\hat\varphi_\gamma}=\frac{\prod_{j=1}^\ell\left|s_{Y_j}\right|^{2\alpha_j}}{\Phi^{\hat q-1-\gamma}\left(\sum_{j=1}^{\tilde\ell}\left|s_j\right|^2\right)^\gamma}
$$
of $$\left(\hat q-1\right)\left(K_X-\sum_{j=1}^\ell a_j Y_j\right)-\sum_{j=1}^\ell a_j Y_j.$$ The condition 
$\alpha_\nu\left(\hat q-1\right)>\gamma\left(\alpha_\nu+\delta_\nu\right)$ for $1\leq\nu\leq\ell^\prime$ is to make sure that the curvature current of the metric $e^{-\hat\varphi_\gamma}$ is a positive current on $X$.  Let $\hat{\mathcal I}_\gamma$ be the multiplier ideal sheaf of the metric $e^{-\hat\varphi_\gamma}$.   Let $W$ denote the subspace of $X$ whose structure sheaf is ${\mathcal O}_X\left/\hat{\mathcal I}_\gamma\right.$  

\medbreak We are going to choose $\gamma$ by imposing more conditions on it, and we are interested in the following two conditions.
\begin{itemize}\item[(a)] $\prod_{j=1}^{\ell^\prime}{\mathcal I}d_{Y_j}$ is contained in $\hat{\mathcal I}_\gamma$ on $X$.
\item[(b)] The restriction of $\hat\rho$ to the subspace ${\mathcal O}_X\left/\left(\hat{\mathcal I}_\gamma+\prod_{j=1}^{\ell^\prime}{\mathcal I}d_{Y_j}\right)\right.$ is not identically zero. In other words, $\hat\rho{\mathcal O}_X$ is not contained in $\hat{\mathcal I}_\gamma+\prod_{j=1}^{\ell^\prime}{\mathcal I}d_{Y_j}$ as ideal sheaves on $X$.
\end{itemize}
In general we cannot hope to be able to choose $\gamma$ such that both Conditions (a) and (b) are satisfied.   If it is possible to choose $\gamma$ satisfying both Conditions (a) and (b), then we can use the general type property of $X$ to slightly modify the metric $e^{-\hat\varphi_\gamma}$ to make its curvature current strictly positive while the multiplier ideal sheaf of the new metric remains the same as $\hat{\mathcal I}_\gamma$.  We then use the vanishing theorem of Kawamata-Viehweg-Nadel to get  the following vanishing of the cohomology
$$
\displaylines{H^1\left(X,\,\hat{{\mathcal I}_\gamma}\left(\hat q\left(K_X-\sum_{j=1}^\ell a_j Y_j\right)\right)\right)\cr=H^1\left(X,\,\hat{{\mathcal I}_\gamma}\left(\left(\left(\hat q-1\right)\left(K_X-\sum_{j=1}^\ell a_j Y_j\right)-\sum_{j=1}^\ell a_j Y_j\right)+K_X\right)\right)=0.}$$
The restriction map
$$\Gamma\left(X,\,\hat q\left(K_X-\sum_{j=1}^\ell a_j Y_j\right)\right)\to
\Gamma\left(X,\,\left({\mathcal O}_X\left/\hat{{\mathcal I}_\gamma}\right.\right)\left(\hat q\left(K_X-\sum_{j=1}^\ell a_j Y_j\right)\right)\right)$$
is surjective.   Let $\rho^\dagger$ be the element of $$\Gamma\left(X,\,\left({\mathcal O}_X\left/\hat{{\mathcal I}_\gamma}\right.\right)\left(\hat q\left(K_X-\sum_{j=1}^\ell a_j Y_j\right)\right)\right)$$ induced by $\hat\rho$.   Condition (a) guarantees that $\rho^\dagger$ is well defined and Condition (b) guarantees that $\rho^\dagger$ is nonzero.  The element $\rho^\dagger$ can be lifted to an element
$$
\rho^\natural\in\Gamma\left(X,\,\hat q\left(K_X-\sum_{j=1}^\ell a_j Y_j\right)\right).
$$
Since $\rho^\dagger$ is nonzero, from Condition (a) it follows that the restriction of $\rho^\natural$ to $\bigcup_{j=1}^{\ell^\prime}Y_j$ is not identically zero. This gives the conclusion that the generic stable vanishing order $\alpha_j$ of $Y_j$ is achieved for some $1\leq j\leq\ell^\prime$ (which actually is a contradiction because we assume that none of the generic stable vanishing order $\alpha_j$ of $Y_j$ is achieved for $1\leq j\leq\ell^\prime$).

\bigbreak\noindent(6.3) {\it Subspaces of Minimal Additional Vanishing.}  We now introduce subspaces of minimal additional vanishing to handle the difficulty of additional vanishing.   We use the same notations as in (6.2).   We blow up $X$ to $\pi: \tilde X\to X$ such that the pullback of $\sum_{j=1}^{\tilde\ell}\left|s_j\right|^2$ to $\tilde X$ is of the form $\prod_{j=1}^{\hat e}\left|s_{E_j}\right|^{2\beta_j}$ for some nonsingular hypersurfaces $E_j$ in $\tilde X$ in normal crossing with each $\beta_j$ being a nonnegative rational number and $K_{\tilde X}-\pi^*K_X$ is a divisor whose support is contained in $\bigcup_{j=1}^{\hat e}E_j$.

\medbreak We now work with $\tilde X$ instead of $X$.  Then $E_j$ is one of $Y_1,\cdots,Y_\ell$ and
$\sum_{j=1}^{\tilde\ell}\left|s_j\right|^2$ is $\prod_{j=1}^{\ell}\left|s_{Y_j}\right|^{2\beta_j}$ with each $\beta_j$ being a nonnegative rational number.  Now for $\nu=\ell^\prime$ the ideal sheaf $\widetilde{{\mathcal I}_{\nu}}$ is precisely equal to ${\mathcal I}_{
\left(q_\nu-1\right)\varphi-q_\nu\sum_{j=1}^\ell\log\left|s_{Y_\nu}\right|^2}$ times the full ideal sheaf ${\mathcal I}d_{Y_\nu}$ of $Y_\nu$.  The difficulty of additional vanishing is no longer a problem for the extension of the section
$$\rho_\nu\in\Gamma\left(Y_\nu,q_\nu\left(K_X-\sum_{j=1}^\ell a_j Y_j\right)\right)$$
to an element of $\Gamma\left(X, q_\nu\left(K_X-\sum_{j=1}^\ell a_j Y_j\right)\right)$ when $\nu=\ell^\prime$.  The subspace $Y_\nu$ with $\nu=\ell^\prime$ is the subspace of ``minimal additional singularity''.

\medbreak This procedure is analogous to the use of the minimum center of log canonical singularity in the arguments for Fujita-conjecture-type problems  [Shokurov 1985, Kawamata 1985, Fujita 1987].    Though ``additional vanishing'' may pose a problem when we insist on doing the extension from a prescribed subspace, yet ``additional vanishing'' can be used to define  another subspace from which extension can be done and in general this new subspace is different from the prescribed subspace.    In the arguments for Fujita-conjecture-type problems the minimum center of log canonical singularities is this new subspace.  In our case this ndew subspace is $Y_\nu$ with $\nu=\ell^\prime$.  

\medbreak  What is attempted in (6.2) is to introduce the subspace $W$ as the subspace of ``minimum additional vanishing'' and {\it require at the same time that $W$ is inside the old hypersurfaces $Y_\nu$ in the old $X$ in (6.2)}.   To go to the subspace of ``minimum additional vanishing'' we have to move out in some direction transversal to the old hypersurfaces $Y_\nu$ in (6.2).   If we insist on the underlying subvariety of $W$ to be inside the old hypersurfaces $Y_\nu$ in (6.2), we would have to consider the situation of $W$ being a subspace of the unreduced space $\left(Y_\nu, {\mathcal O}_X\left/\left({\mathcal I}d_{Y_\nu}\right)^b\right.\right)$ for some integer $b>1$, though $W$ may be a subspace inside {\it another reduced hypersurface}  $Y^\prime$.  Only after blowing up, we can make $W$ inside the new hypersurface $Y_\nu$ and, as a matter of fact, even precisely equal to one of them.  That is precisely what we are doing here in blowing up $X$ to $\tilde X$.

\medbreak One problem with this kind of blow-up argument to get the subspace of ``minimum additional vanishing'' is that we have to worry about the process not terminating before we arrive at any of our original hypersurfaces of the second case of the dichotomy.   Such a termination is obtained by using the technique of discrepancy subspaces in \S4 and the technique of the construction of metrics of additional singularity on hypersurface of first case of dichotomy after fixed ample twisting in \S5.

\bigbreak\noindent(6.4) {\it  Shokurov's Technique of Comparing the Theorem of Hirzebruch-Riemann-Roch for a Line Bundle and for Its Flat-Twisting.}  We now consider Item (i) of (6.2).   It has the same difficulty of ``additional vanishing'' as in Item (ii) of (6.2).  It can simply be dealt with in the same way as in Item (ii) of (6.3) by using the technique of subspaces of ``minimum additional vanishing'' presented in (6.3).  For $1\leq\nu\leq\ell^\prime$ we take the metric
$$
e^{-\tilde\varphi_\nu}=\frac{\left(\prod_{j=1}^\ell \left|s_{Y_j}\right|^{2a_j}\right)^{q_\nu}}{\Phi^{q_\nu-1-\frac{1+\alpha_\nu}{\delta_\nu}}\left(\sum_{j=1}^{\tilde\ell}\left|s_j\right|^2\right)^{\frac{1+\alpha_\nu}{\delta_\nu}}}
$$
of $$\left(q_\nu-1\right)\left(K_X-\sum_{j=1}^\ell a_j Y_j\right)-\sum_{j=1}^\ell a_j Y_j$$ on $X$, which was introduced in (6.2) and whose multiplier ideal sheaf 
$\widetilde{{\mathcal I}_\nu}$ is equal to ${\mathcal I}_{
\left(q_\nu-1\right)\varphi-q_\nu\sum_{j=1}^\ell\log\left|s_{Y_\nu}\right|^2}\,{\mathcal I}d_{Y_\nu}$  when $\nu=\ell^\prime$.   We can now define the metric
$$
e^{-\varphi_\nu^\sharp}=\left.\frac{\left(\prod_{j=1}^\ell \left|s_{Y_j}\right|^{2a_j}\right)^{q_\nu}\left|s_{Y_\nu}\right|^2}{\Phi^{q_\nu-1-\frac{1+\alpha_\nu}{\delta_\nu}}\left(\sum_{j=1}^{\tilde\ell}\left|s_j\right|^2\right)^{\frac{1+\alpha_\nu}{\delta_\nu}}}\right|_{Y_\nu}
$$
of $$\left(q_\nu-1\right)\left(K_X-\sum_{j=1}^\ell a_j Y_j\right)-\sum_{j=1}^\ell a_j Y_j-Y_\nu$$ on $Y_\nu$ when $\nu=\ell^\prime$, whose multiplier ideal sheaf 
${{\mathcal I}_\nu}^\sharp$ on $Y_\nu$  is equal to ${\mathcal I}_{\left.\left(
\left(q_\nu-1\right)\varphi-q_\nu\sum_{j=1}^\ell\log\left|s_{Y_\nu}\right|^2\right)\right|_{Y_\nu}}$  when $\nu=\ell^\prime$.   By using the general property of $X$ to slightly change the metric $e^{-\varphi_\nu^\sharp}$ on $Y_\nu$, we can assume (by keeping the same symbol) that the curvature current of the metric  $e^{-\varphi_\nu^\sharp}$ on $Y_\nu$ is strictly positive for $\nu=\ell^\prime$.  By applying the vanishing theorem of Kawamata-Viehweg-Nadel to $Y_\nu$, we conclude for $\nu=\ell^\prime$ that
$$
H^\lambda\left(Y_\nu, \,{{\mathcal I}_\nu}^\sharp\left(q_\nu\left(K_X-\sum_{j=1}^\ell a_j Y_j\right)\right)\right)=0\quad{\rm for\ \ }\lambda\geq 1
$$
because $K_{Y_\nu}=K_X+Y_\nu$.  Moreover, for any flat bundle $F$ over $Y_\nu$ for $\nu=\ell^\prime$ we also have
$$
H^\lambda\left(Y_\nu, \,{{\mathcal I}_\nu}^\sharp\left(q_\nu\left(K_X-\sum_{j=1}^\ell a_j Y_j\right)\right)+F\right)=0\quad{\rm for\ \ }\lambda\geq 1.
$$
Since for $\nu=\ell^\prime$ the modified restriction to $Y_\nu$ of the curvature current of the metric $\frac{1}{\Phi}$ of $K_X$ is equal to $\sum_{j=1}^{J_\nu}\gamma_{\nu,j}\left[V_{\nu,j}\right]$, it follows that the section
$\prod_{j=1}^{J_\nu}\left(s_{V_{\nu,j}}\right)^{\gamma_{\nu,j}q_\nu}$ is a holomorphic section of $\left.q_\nu\left(K_X-\alpha_\nu Y_\nu\right)\right|_{Y_\nu}+F$  over $Y_\nu$ for some flat line bundle $F$ over $Y_\nu$.   Clearly $\prod_{j=1}^{J_\nu}\left(s_{V_{\nu,j}}\right)^{\gamma_{\nu,j}q_\nu}$ can be naturally regarded as a nonzero element of $$\Gamma\left(Y_\nu, \,{{\mathcal I}_\nu}^\sharp\left(q_\nu\left(K_X-\sum_{j=1}^\ell a_j Y_j\right)+F\right)\right)$$  for $\nu=\ell^\prime$.  By comparing the results of the application of the theorem of Hirzebruch-Riemann-Roch respectively to ${{\mathcal I}_\nu}^\sharp\left(q_\nu\left(K_X-\sum_{j=1}^\ell a_j Y_j\right)+F\right)$ and to ${{\mathcal I}_\nu}^\sharp\left(q_\nu\left(K_X-\sum_{j=1}^\ell a_j Y_j\right)\right)$ on $Y_\nu$ for $\nu=\ell^\prime$, we conclude that there is a nonzero element
$$
\rho_\nu\in\Gamma\left(Y_\nu, \,{{\mathcal I}_\nu}^\sharp\left(q_\nu\left(K_X-\sum_{j=1}^\ell a_j Y_j\right)\right)\right)$$  for $\nu=\ell^\prime$.

\bigbreak\noindent\S7. {\it Big Sum of a Line Bundle and the Canonical Line Bundle.}  For finite generation of the canonical ring when $X$ is not of general type, the case of adding a  line bundle $L$ over $X$ to the canonical line bundle $K_X$ with the sum being big can be essentially considered as the case of just the canonical line bundle over the total bundle space of the dual $L^*$ of $L$ (when only pluricanonical sections whose coefficients with respect to the fiber coordinates of $L^*$ are constant along the fibers of $L^*$ are considered).  The reason is as follows.

\medbreak Let $z^1_j,\cdots,z^n_j$ be the local holomorphic coordinate system on an open subset $U_j$ of $X$ and $g_{jk}$ be the transition function for the line bundle $L^*$ from the fiber coordinate $w_k$ on $U_k$ to the fiber coordinate $w_j$ on $U_j$ so that $w_j=g_{jk}w_k$. Then $dw_j=\left(dg_{jk}\right)w_k+g_{jk}dw_k$ and
$$
dz^1_j\wedge\cdots\wedge dz^n_j\wedge dw_j=\frac{\partial\left(z^1_j,\cdots,z^n_j\right)}
{\partial\left(z^1_k,\cdots,z^n_k\right)}\,
dz^1_k\wedge\cdots\wedge dz^n_k\wedge g_{jk}dw_k
$$
over $U_j\cap U_k$.
For any element $s\in\Gamma\left(X,m\left(L+K_X\right)\right)$ represented by $\left\{s_j\right\}_j$ with $s_j$ being a holomorphic function on $U_j$, we have 
$$
s_j\left(dz^1_j\wedge\cdots\wedge dz^n_j\wedge dw_j\right)^{\otimes m}=
s_k\left(dz^1_k\wedge\cdots\wedge dz^n_k\wedge dw_k\right)^{\otimes m}
$$
over $U_j\cap U_k$.  This means that when we consider only elements of $\Gamma\left(L^*,mK_{L^*}\right)$ whose coefficients with respect to the fiber coordinates of $L^*$ are constant along the fibers of $L^*$, we can assume that the canonical line bundle $K_{L^*}$ of $L^*$ to be big and $L^*$ is of general type (in the sense that only elements whose coefficients with respect to the fiber coordinates of $L^*$ are constant along the fibers of $L^*$ are being considered).  Then we can modify the analytic proof of the finite generation of the canonical ring for a compact complex algebraic manifold of general type to give us the analytic proof of the finite generation of the canonical ring over $L^*$ when only pluricanonical sections over $L^*$ whose coefficients with respect to the fiber coordinates of $L^*$ are constant along the fibers of $L^*$ are being considered.

\bigbreak\noindent{\bf APPENDIX: Multiplier Ideal Sheaves of Kohn and Nadel as Defined by Crucial Estimates}

\bigbreak As mentioned at the beginning of the Introduction and in (4.1) the definition of discrepancy  subspaces is motivated by the original philosophy of formulating multiplier ideal sheaves from the most crucial estimates when multiplier ideal sheaves were first introduced by Kohn as measurements of failure of estimates in partial differential equations [Kohn 1979] and introduced by Nadel as destabilizing sheaves [Nadel 1990].  In this Appendix we examine the definitions of the original multiplier ideal sheaves of Kohn and Nadel and, especially, recast Nadel's definition in the context of formulation in terms of the most crucial estimates.

\bigbreak\noindent(A.1) {\it Kohn's Subelliptic Multipliers
for the Complex Neumann Problem.}  The setting of Kohn's multiplier ideal sheaf is a bounded domain $\Omega$ in ${\bf C}^n$ with smooth
weakly pseudoconvex boundary defined by $r<0$ with $dr$ being nowhere
zero on the boundary $\partial\Omega$ of $\Omega$. Here weakly pseudoconvex boundary means that
$\sqrt{-1}\,\partial\bar\partial r|_{T^{(1,0)}_{\partial\Omega}}\geq 0$.  The problem is to study the following regularity question: given a smooth $(0,1)$-form $f$ on $\bar\Omega$ with $\bar\partial f=0$, whether the solution of
$\bar\partial u=f$ on $\Omega$ with $u$ perpendicular to all holomorphic functions on
$\Omega$ is smooth on $\bar\Omega$.

\medbreak A sufficient condition for regularity is the following subelliptic estimate at every boundary point.  For $P\in\partial \Omega$ there exist some
open neighborhood $U$ of $P$ in ${\bf C}^n$ and positive numbers
$\epsilon$ and $C$ satisfying
$$
\||g|\|_\epsilon^2\leq C\left(\|\bar\partial g\|^2+\|\bar\partial^*
g\|^2+\|g\|^2\right)\leqno{({\rm A}.1.1)}
$$
for every $(0,1)$-form $g$ supported on $U\cap\bar\Omega$ which is in the
domain of $\bar\partial$ and $\bar\partial^*$.  Here $\||\cdot|\|_\epsilon$ is the $L^2$ norm on $\Omega$
involving derivatives up to order $\epsilon$ in the boundary
tangential directions of $\Omega$, $\|\cdot\|$ is the usual $L^2$ norm on $\Omega$ without involving any
derivatives, and $\bar\partial^*$ is the actual adjoint of
$\bar\partial$ with respect to $\|\cdot\|$.

\medbreak The reason why some positive $\varepsilon$ is needed is that in applying a differential operator $D$ to both sides of $\bar\partial u=f$ to get estimates of the Sobolev norm of $u$ up to a certain order of derivatives in terms of that of $f$, an error term from the commutator of the differential operator $D$ and $\bar\partial$ occurs, which needs to be absorbed and one way to do the absorption is to use an estimate involving a Sobolev norm with derivative higher by some positive number $\varepsilon$.  This stronger Sobolev norm is used also to absorb the error term from partitions of unity or cut-off functions.

\medbreak The reason why only the tangential Sobolev norm $\||\cdot|\|_\epsilon$ is used is that we need to preserve the condition that $(0,1)$-form $g$ belongs to the domain of $\bar\partial^*$ (which means vanishing complex-normal component at boundary points) by using only differentiation along the boundary tangential directions.  The missing estimate in the real-normal direction can be obtained from the complex-normal component of the equation $\bar\partial u=f$.

\medbreak The theory of multiplier ideal sheaves introduces multipliers into the most crucial estimate, which in this case is the subelliptic estimate (A.1.1).   A {\it subelliptic multiplier} $F$ is a smooth function germ of ${\mathbb C}^n$ at $P$ such that the following subellitpic
estimate of order $\varepsilon_F$ holds for any test $(0,1)$-form $g$
after replacing it by its product with $F$.
$$
\||Fg|\|_{\epsilon_{{}_F}}^2\leq C_{{}_F}\left(\|\bar\partial
g\|^2+\|\bar\partial^* g\|^2+\|g\|^2\right)\leqno{({\rm A}.1.2)}
$$
for every test $(0,1)$-form $g$ described above.  The {\it multiplier ideal} $I_P$ at the boundary point $P$ is the ideal of all such subelliptic multipliers $F$.    Kohn's multiplier ideal sheaf is the sheaf of such ideals $I_P$. 

\medbreak Since each test $(0,1)$-form $g$ presents one inequality (A.1.2), Kohn's multiplier ideal sheaf is actually defined by a family of inequalities (A.1.2), parametrized by the set of all test $(0,1)$-forms $g$.   Kohn's multiplier ideal sheaf is therefore a {\it dynamic} multiplier ideal sheaf.

\medbreak In Kohn's case it is clear that the definition (A.1.2) of the multiplier ideal sheaf is formulated from the most crucial estimate (A.1.1).  However, the situation is by no means clear in Nadel's original formulation of his multiplier ideal sheaves.   Here we are going to recast Nadel's multiplier ideal sheaves in the context of formulation in terms of the most crucial estimates.

\bigbreak\noindent(A.2) {\it Nadel's Multiplier Ideal Sheaves.}  The setting of Nadel's multiplier ideal sheaves is a compact complex manifold $X$ of complex dimension $n$ with an ample anticanonical line bundle $-K_X$ of $X$ .  Let $g_{i\bar j}$
be a K\"ahler metric of $X$ in the anticanonical class of $X$.  Let
$$
R_{i\bar j}=-\partial_i\partial_{\bar j}\det\left(g_{i\bar j}\right)_{1\leq i,j\leq n}
$$ be the Ricci curvature of $g_{i\bar j}$.  There is a smooth positive function $F$ on $X$ such that
$$
R_{i\bar j}-g_{i\bar j}=\partial_i\partial_{\bar j}\log F.
$$
We consider the
complex Monge-Amp\`ere equation
$$
\det\left(g_{i\bar j}+\partial_i\partial_{\bar j}\varphi\right)_{1\leq i,j\leq n}
=e^{-\varphi}F\det\left(g_{i\bar j}\right)_{1\leq i,j\leq n},\leqno{({\rm A}.2.1)}
$$
formulated by Calabi [Calabi 1954a, Calabi 1954b, Calabi 1955] for the construction of a K\"ahler-Einstein metric of $X$.  If the equation (A.2.1) is solved, by taking $\partial\bar\partial\log$ of both sides of (A.2.1), we get
$$
-R^\prime_{i\bar j}=-\left(g^\prime_{i\bar j}-g_{i\bar j}\right)+\left(R_{i\bar j}-g_{i\bar j}\right)-R_{i\bar j}=-g^\prime_{i\bar j},
$$
(where $g^\prime_{i\bar j}=g_{i\bar j}+\partial_i\partial_{\bar j}\varphi$ and $R^\prime_{i\bar j}$ is the Ricci curvature of the K\"ahler metric $g^\prime_{i\bar j}$) and conclude that $g^\prime_{i\bar j}$ is a K\"ahler-Einstein metric of $X$.  Continuity method is applied to solve the equation (A.2.1) by considering the solution of
$$
\det\left(g_{i\bar j}+\partial_i\partial_{\bar j}\varphi_t\right)_{1\leq i,j\leq n}
=e^{-t\varphi_t}F\det\left(g_{i\bar j}\right)_{1\leq i,j\leq n},\leqno{({\rm A}.2.2)_t}
$$
for $0\leq t\leq 1$, starting with $t=0$ by using [Yau 1978, p.363, Theorem 1].

\medbreak The openness part of the continuity method is clear from the usual elliptic estimates and the implicit function theorem.  Nadel's multiplier ideal sheaf arises from the closedness part of the continuity method in the following way.  Suppose for some $0<t_*\leq 1$ we have a sequence $\varphi_{t_\nu}$ which satisfies $({\rm A}.2.2)_{t_\nu}$ with $t_\nu\to t_*$ monotonically strictly increasing as $\nu\to\infty$.

\medbreak \medbreak Since the first Chern class of $-K_X$, which (up to a normalizing universal constant) is represented by
$$\sum_{i,j=1}^n\left(g_{i\bar j}+\partial_i\partial_{\bar j}\varphi_t\right)\left(\frac{\sqrt{-1}}{2}dz_i\wedge d\overline{z_j}\right),\leqno{({\rm A}.2.3)_t}
$$
is independent of $t<t_*$, the $(1,1)$-form $({\rm A}.2.3)_t$ would converge weakly when $t$ goes through an appropriate sequence $t_\nu$ to $t_*$.  Let $\widehat{\varphi_t}$ be the average of $\varphi_t$ over $X$ with respect to the K\"ahler metric $g_{i\bar j}$.  Since the Green's operator for the Laplacian, with respect to the K\"ahler metric $g_{i\bar j}$, is a compact operator from the space of bounded measures on $X$ to the space of $L^1$ functions on $X$, we conclude that $\varphi_{t_\nu}-\widehat{\varphi_{t_\nu}}$ converges to some function in the $L^1$ norm for some subsequence $t_\nu$ of $t\to t_*$.

\medbreak The second-order and third-order estimates used to obtain [Yau1978, p.363, Theorem 1] work also for applying the continuity method to solve $({\rm A}.2.2)_t$ for $0\leq t\leq 1$.  Alternatively the H\"older estimate for the second-order derivatives can be used instead of the third-order estimates (see {\it e.g.,} [Siu 1987, Chapter 2, \S3 and \S4]).

\medbreak The obstacle in the closedness part $t\to t_*$ of the continuity method for solving $({\rm A}.2.2)_t$ occurs when $\widehat{\varphi_{t_\nu}}\to\infty$ as $\nu\to\infty$.  After multiplying $({\rm A}.2.2)_{t_\nu}$ by $e^{t_\nu\widehat{\varphi_{t_\nu}}}$ to get
$$
e^{t_\nu\widehat{\varphi_{t_\nu}}}\det\left(g_{i\bar j}+\partial_i\partial_{\bar j}\varphi_{t_\nu}\right)_{1\leq i,j\leq n}
=e^{-t_\nu\left(\varphi_{t_\nu}-\widehat{\varphi_{t_\nu}}\right)}F\det\left(g_{i\bar j}\right)_{1\leq i,j\leq n}
$$
and integrating over $X$ and taking limit as $\nu\to\infty$, we get
$$
\lim_{\nu\to\infty}\int_X e^{-t_\nu\left(\varphi_{t_\nu}-\widehat{\varphi_{t_\nu}}\right)}=\infty\leqno{({\rm A}.2.4)}
$$
when $\widehat{\varphi_{t_\nu}}\to\infty$ as $\nu\to\infty$,
because
$$
\displaylines{\int_X\det\left(g_{i\bar j}+\partial_i\partial_{\bar j}\varphi_{t_\nu}\right)_{1\leq i,j\leq n}\prod_{j=1}^n\left(\frac{\sqrt{-1}}{2}dz_j\wedge d\overline{z_j}\right)\cr=
\int_X\det\left(g_{i\bar j}\right)_{1\leq i,j\leq n}\prod_{j=1}^n\left(\frac{\sqrt{-1}}{2}dz_j\wedge d\overline{z_j}\right)=
\left(-K_X\right)^n\cr}
$$
which is independent of $t$.

\medbreak We now know that the crucial estimate in Nadel's setting is
$$
\lim_{\nu\to\infty}\int_X e^{-t_\nu\left(\varphi_{t_\nu}-\widehat{\varphi_{t_\nu}}\right)}<\infty.
$$
Since the multiplier ideal sheaf is introduced to make the crucial estimate hold after using a multiplier, we introduce the multiplier ideal sheaf ${\mathcal I}$ in Nadel's setting as consisting of all holomorphic function germs $f$ on $X$ such that
$$
\limsup_{\nu\to\infty}\int_U \left|f\right|^2e^{-t_\nu\left(\varphi_{t_\nu}-\widehat{\varphi_{t_\nu}}\right)}<\infty,
$$
where $U$ is an open neighborhood of the point of $X$ at which $f$ is a germ.  This multiplier ideal sheaf ${\mathcal I}$ in the sense of Nadel is defined by using a sequence of functions $\varphi_{t_\nu}-\widehat{\varphi_{t_\nu}}$ as $\nu\to\infty$ and is therefore a {\it dynamic} multiplier ideal sheaf.

\bigbreak\noindent{\it References}

\medbreak\noindent[Birkan-Cascini-Hacon-McKernan 2006] C. Birkar, P.
Cascini, C. Hacon, and J. McKernan, Existence of minimal models for
varieties of log general type, arXiv:math/0610203.

\medbreak\noindent[Calabi 1954a] E. Calabi, The variation of K\"ahler metrics I: The
structure of the space; II: A minimum problem, {\it Amer. Math.
Soc. Bull.} \textbf{60} (1954), Abstract  Nos. 293--294, p.168.

\medbreak\noindent[Calabi 1954b] E. Calabi, The space of K\"ahler metrics, {\it
Proc. Internat. Congress Math}.  Amsterdam, 1954, Vol. 2,
pp.206--207.

\medbreak\noindent[Calabi 1955] E. Calabi, On K\"ahler manifolds with vanishing
canonical class, {\it Algebraic Geometry and Topology, A Symposium
in Honor of S. Lefschetz}, Princeton Univ. Press, Princeton, 1955,
pp.78--89.

\medbreak\noindent[Demailly 1992] J.-P. Demailly, Regularization of
closed positive currents and Intersection Theory, {\it J. Alg.
Geom.} {\bf 1} (1992) 361-409.

\medbreak\noindent
[Fujita 1987] T. Fujita, On polarized manifolds whose adjoint bundles are not
semipositive, Proceedings of the 1985 Sendai Conference on Algebraic Geometry,
{\it Advanced Studies in Pure Mathematics} \textbf{10} (1987), 167-178.

\medbreak\noindent[Hardy-Wright 1960] G. H. Hardy and E. M. Wright,
{\it An Introduction to the Theory of Numbers}, 4th ed., Oxford
University Press 1960.

\medbreak\noindent[Kawamata 1982] Y. Kawamata, A generalization of
Kodaira-Ramanujam's vanishing theorem. {\it Math. Ann.} {\bf 261}
(1982), 43-46.

\medbreak\noindent[Kawamata 1985] Y. Kawamata, Pluricanonical
systems on minimal algebraic varieties. {\it Invent. Math.} {\bf 79}
(1985), 567--588.

\medbreak\noindent [Kohn 1979] Joseph J. Kohn, Subellipticity of the
$\bar \partial $-Neumann problem on pseudo-convex domains:
sufficient conditions. {\it Acta Math.} \textbf{142} (1979), 79--122.

\medbreak\noindent [Lelong 1968] Pierre Lelong,
{\it Fonctions plurisousharmoniques et formes diff\'erentielles positives}. Gordon \& Breach, Paris-London-New York (Distributed by Dunod \'editeur, Paris) 1968.

\medbreak\noindent[Nadel 1990] A. Nadel, Multiplier ideal sheaves
and K\"ahler-Einstein metrics of positive scalar curvature. {\it
Ann. of Math.} {\bf 132} (1990), 549-596.

\medbreak\noindent[Ohsawa-Takegoshi 1987] Ohsawa, T., Takegoshi, K.:
On the extension of $L^2$ holomorphic functions, {\it Math.
Zeitschr.} {\bf 195}, 197-204 (1987).

\medbreak\noindent[Paun 2005] M. Paun, Siu's invariance of
plurigenera: a one-tower proof, {\it J. Differential Geom.} {\bf 76}
(2007), no. 3, 485--493.

\medbreak\noindent[Paun 2005] M. Paun,  Relative critical exponents, non-vanishing and metrics with minimal singularities,  arXiv:0807.3109.

\medbreak\noindent[Shokurov 1985] V.~V. Shokurov, A nonvanishing
theorem. {\it Izv. Akad. Nauk SSSR Ser. Mat.} {\bf 49} (1985),
635--651.

\medbreak\noindent[Siu 1974] Y.-T. Siu, Analyticity of sets
associated to Lelong numbers and the extension of closed positive
currents. {\it Invent. Math.} {\bf 27} (1974), 53-156.

\medbreak\noindent[Siu1987]
Yum-Tong Siu,
{\it Lectures on Hermitian-Einstein metrics for stable bundles and K\"ahler-Einstein metrics.}
DMV Seminar, 8. Birkh\"auser Verlag, Basel, 1987.

\medbreak\noindent[Siu 1996] Siu, Y.-T., The Fujita conjecture and
the extension theorem of Ohsawa-Takegoshi, in {\it Geometric Complex
Analysis} ed. Junjiro Noguchi {\it et al}, World Scientific:
Singapore, New Jersey, London, Hong Kong 1996, pp. 577-592.

\medbreak\noindent[Siu 1998]  Y.-T. Siu, Invariance of plurigenera,
{\it Invent. Math.} {\bf 134} (1998), 661-673.

\medbreak\noindent[Siu 2002]  Y.-T. Siu, Extension of twisted
pluricanonical sections with plurisubharmonic weight and invariance
of semipositively twisted plurigenera for manifolds not necessarily
of general type. In: {\it Complex Geometry: Collection of Papers
Dedicated to Professor Hans Grauert}, Springer-Verlag 2002,
pp.223-277.

\medbreak\noindent[Siu 2003] Y.-T. Siu, Invariance of plurigenera
and torsion-freeness of direct image sheaves of pluricanonical
bundles, In: {\it Finite or Infinite Dimensional Complex Analysis
and Applications} (Proceedings of the 9th International Conference
on Finite or Infinite Dimensional Complex Analysis and Applications,
Hanoi, 2001), edited by Le Hung Son, W. Tutschke, C.C. Yang, Kluwer
Academic Publishers 2003, pp.45-83.

\medbreak\noindent[Siu 2005]  Y.-T. Siu, Multiplier ideal sheaves in
complex and algebraic geometry, {\it Science in China,
Ser.A:Math.}{\bf 48} (2005), 1-31 (arXiv:math.AG/0504259).

\medbreak\noindent[Siu 2006] Y.-T. Siu, A general non-vanishing
theorem and an analytic proof of the finite generation of the
canonical ring, arXiv:math/0610740.

\medbreak\noindent[Siu 2007] Y.-T. Siu, Additional explanatory notes
on the analytic proof of the finite generation of the canonical
ring, arXiv:0704.1940.

\medbreak\noindent[Siu 2008] Y.-T. Siu,
Finite generation of canonical ring by analytic method (arXiv:0803.2454),
{\it J. Sci. China} \textbf{51} (2008), 481-502.

\medbreak\noindent[Skoda 1972] H. Skoda, Application des techniques
$L^2$ \`a la th\'eorie des id\'eaux d'une alg\`ebre de fonctions
holomorphes avec poids. {\it Ann. Sci. \'Ecole Norm. Sup.} {\bf 5}
(1972), 545-579.

\medbreak\noindent[Viehweg 1982] E. Viehweg, Vanishing theorems.
{\it J. Reine Angew. Math.} {\bf 335} (1982), 1-8.

\medbreak\noindent[Yau 1978]
Shing-Tung Yau,
On the Ricci curvature of a compact K\"ahler manifold and the complex Monge-Amp\`ere equation. I.
{\it Comm. Pure Appl. Math.} \textbf{31} (1978), 339--411.

\bigbreak\noindent{\it Author's mailing address}: Department of
Mathematics, Harvard University, Cambridge, MA 02138, U.S.A.

\medbreak\noindent {\it Author's e-mail address}:
siu@math.harvard.edu

\end{document}